\input amstex
\documentstyle{amsppt}\magnification=\magstep1
\magnification=1200
\NoRunningHeads
\NoBlackBoxes
\baselineskip=12pt
\vsize=8.5truein
\hsize=6.0truein
\hoffset=0.3truein
\voffset=0.2truein
\loadmsam
\loadbold
\loadmsbm
\loadeufm \loadeusm \UseAMSsymbols
 \NoRunningHeads
 \nologo
\document
\define\({\left(}
\define\){\right)} \define\[{\left[} \define\]{\right]}
\define\db{\bar\partial} \define\p{\parallel}
\define\aint{\rlap{\kern2pt\vrule height2.5pt width9.5pt depth-2.1pt}
\int_{B_{\rho}}}
   \nologo

\topmatter
\title Estimates for the $\db$-Neumann problem     and nonexistence of
 Levi-flat hypersurfaces  in $\Bbb CP^n$
\endtitle
\author   Jianguo Cao,     Mei-Chi Shaw  and Lihe Wang  \endauthor
\author   Jianguo Cao\footnote{ Partially
supported by NSF grant DMS-0102552.
J. Cao is   grateful to the University of Michigan 
 for the hospitality.  Research at MSRI is supported in part by NSF grant DMS-9810361.  \hfill{$\,$}  },     Mei-Chi
Shaw\footnote{   Partially
supported by NSF grant DMS 01-00492.   \hfill{$\,$}   }
 and Lihe Wang\footnote{ Partially
supported by NSF grant DMS 01-00679.   \hfill{$\,$} }
\endauthor
\address Jianguo Cao,
Department of Mathematics, University of Notre Dame, Notre Dame, IN  46556
USA  \endaddress
\email cao.7$\@$nd.edu\endemail
\address Mei-Chi Shaw,
Department of Mathematics, University of Notre Dame, Notre Dame, IN  46556
USA  \endaddress
\email mei-chi.shaw.1\@nd.edu \endemail
\address Lihe Wang,
 Department of Mathematics, University of Iowa,
      Iowa City, IA 52242 USA
\endaddress
\email lwang\@math.uiowa.edu \endemail

\abstract
Let $\Omega$ be a {\it  pseudoconvex } domain with $C^2$-smooth boundary
in $\Bbb CP^n$. We prove that the  $\db$-Neumann operator $N$ exists for
$(p,q)$-forms  on
$\Omega$.
Furthermore, there exists a $t_0>0$ such that  the operators
 $N$, $\db^*N$, $\db N$ and the Bergman projection
  are regular in the Sobolev space  $W^t (\bar{\Omega}) $ for   $t<t_0$.

 The boundary estimates above have applications in complex geometry.
We use the estimates to prove  the nonexistence of $C^{2, \alpha}$
real
 Levi-flat hypersurfaces  in $\Bbb  CP^n$. We also  show that there exist no
non-zero $L^2$-holomorphic $(p, 0)$-forms on any {\it pseudoconcave}
domain in $\Bbb CP^n$ with $p > 0$.

\endabstract

\endtopmatter

\heading   {\bf Introduction} \endheading

One of the main results in this paper is the following:
\proclaim{Theorem 1} There exists no $C^{2,\alpha}$ real Levi-flat
hypersurface $M^{2n-1}$ in $\Bbb CP^n$, where
$  0<\alpha<1$ and
$n\ge 2$.
\endproclaim
Theorem 1 is inspired by the
recent papers of Siu [Siu2,3] who proved that there exists no
$C^8$ Levi-flat hypersurface in $\Bbb C P^n$, $n\ge 2$.  The required
smoothness has been reduced to $C^4$ by  Iordan [Io].
   Nonexistence of real analytic
Levi-flat hypersurfaces in $\Bbb CP^n$ was
obtained in Lins Neto [LNe] for $n\ge 3$ and Ohsawa [Oh] for $n=2$. Our
proof of Theorem 1 follows  arguments along the lines  of  \cite{Siu2,3}
who reduced the proof of Theorem 1 to the regularity of the tangential
Cauchy-Riemann equations on $M$.

We    derive
 new  boundary regularity results
            for $\db$-equation and  the $\db$-Neumann problem  for domains 
in $\Bbb CP^n$. Our  Theorem 2 stated below 
can be applied to {\it not  only}  domains $\Omega$ with Levi-flat boundaries $M $  {\it but also   all}
other weakly
 pseudoconvex $C^2$ domains in $\Bbb CP^n$, when compared with the earlier work mentioned above.  The
boundary  regularity of the solutions  is  interesting  itself.

To prove Theorem 1,  Siu [Siu2,3] made the following  observation:   If there
exists a $C^2$-smooth {\it real } Levi-flat hypersurface $M = M^{2n-1}$ in
$\Bbb CP^n$ and if
$\Bbb CP^n$ has the standard Fubini-Study metric, then  the curvature form
$i\tilde{\Theta}^N$ of its complex normal
line bundle
$\Cal N_{M, \Bbb CP^n}$ is {\it  strictly positive definite} on
$ T^{(1, 0)}(M)
\oplus T^{(0, 1)}(M)
$ by the Cartan-Chern-Gauss structure equation.

 Furthermore,
  one has $\tilde{\Theta}^N= d\theta$ is an exact form,  where
 $\theta$ is the connection form of the complex normal line bundle $\Cal
N_{M, \Bbb CP^n}$.
If $M$ is a Levi-flat hypersurface, it    is locally
  foliated by complex manifolds of complex dimension $n-1$.
  It is known that  the restriction of $\tilde{\Theta}^N$ to each complex
leaf of $M$,  $\Theta_b =
\tilde{\Theta}^N|_{[T(M)]_{\Bbb R}
 \cap J [T(M)]_{\Bbb R} }$ is a  $(1, 1)$-form, where $J$ is the complex
structure of $\Bbb CP^n$ (see Proposition 1.2
below).

On the other hand, if one could
show that   $\Theta_b  $  is  $\partial_b\db_b$-exact for some
continuous real-valued
function
$h$ (i.e., $\Theta_b=  \partial_b\db_b h$ on $M$), then
$i\Theta_b$ is {\it non-positive  } at the maximum point of
h on $M$ (cf. [Ca]). This  contradicts  the fact that $i\tilde{\Theta}^N$
is {\it  strictly positive definite} on
$ T^{(1, 0)}(M)
\oplus T^{(0, 1)}(M)
$,
 and hence Theorem 1 would follow
immediately.

Therefore,
 the proof of Theorem 1 is reduced to
a problem of finding a continuous function $u $ satisfying
$$
i\partial_b\db_b u = f_b = f|_{  [T(M)]_{\Bbb R}
\cap J [T(M)]_{\Bbb R}  } \quad \text{ on } M,  \tag0.1$$
  under the condition that $f = d\theta$ is an exact real-valued
(1, 1)-form when restricted to $T^{(1,0)}(M) \oplus T^{(0,1)}(M)$. Equation
(0.1)    corresponds to the
classical $\partial\db$-Lelong equation
$$
  i \partial\db \tilde u = \tilde f \quad \text{ in } \Bbb CP^n \tag0.2
$$
   where $\tilde f = d \tilde \theta$ is an exact real-valued (1, 1)-form.
Using the fact that $\tilde f$ is an exact real-valued (1,1)-form
to solve equation
$(0.2)$,  it suffices to solve
$$\db\tilde u=\tilde \theta^{(0, 1)} \quad \text{in } \Bbb CP^n,\tag 0.3$$
where   $\tilde \theta^{(0, 1)}   $ is the
$(0,1)$ part in $\tilde\theta$, which is $\db$-closed in $\Bbb CP^n$.
  By the  Hodge theory,
there is no nontrivial harmonic $(0,1)$-form  in $\Bbb CP^n $.
It follows that any
$\db$-closed (0,1)-form
$ \tilde \theta^{(0, 1)}  $ on $\Bbb CP^n$
must be $\db$-exact, and hence equation (0.3) can be solved easily (e.g., cf.
\cite{Zh}). For the proof of Theorem 1 and equation (0.1), one can
similarly  deduce that  it   suffices
to solve
$$\db_b u= \theta_b \quad \text{in } M,\tag 0.4$$
where $ \theta_b$ is a $(0,1)$-form in $M$ satisfying some compatibility
condition.
In \cite{Siu2,3},    existence and regularity
of equation (0.4) on a Levi-flat boundary $M$ in $\Bbb CP^n$ is studied. In
this paper, we study the more general
situation when $M$ is the $C^2$ boundary of any pseudoconvex domain $\Omega$
in $\Bbb CP^n$.

When $n>2$, the  compatibility condition for equation (0.4)  to be solvable
is that
$$\db_b \theta_b = 0\quad \text{on }M.\tag 0.5$$  When $n=2$, (0.5) is
satisfied trivially and the
compatibility condition for equation (0.4)   is substituted  by  the
following moment
condition:
$$\int_M \theta_b\wedge \Psi =0, \tag 0.6$$
where $\Psi$ is  any $\db_b$-closed  (2,0)-form  on
$M  $  (see (9.2.12a) of [CS, p216].) To show that such a
compatibility condition holds for our case with $n = 2$,
we also derived several
Liouville type theorems which are of independent interest.
In particular, we show that  on any
 {\it pseudoconcave} domain $\Omega $ with $C^2$ boundary in $\Bbb CP^n$,
 there exist no non-zero $L^2$-integrable holomorphic $(p, 0)$-forms on
 $\Omega $ with $p > 0$,
(cf. Proposition 4.5 below).

If one can
 extend a  $(0,1)$-form $\theta_b$ satisfying (0.5) or (0.6) from the
boundary $M = b\Omega$ to a   $\db$-closed form $
 \tilde  \theta $ on the domain $\Omega$,  then   (0.4) can be solved by
restricting the
solution $\tilde u$  of the $\db$-equation (0.3) to $M = b\Omega$.
 To study  the $\db$-closed extension from $M$ to $\Omega$, one can
formulate it as a $\db$-Cauchy problem of  finding
solutions to the $\db$-equation (0.3) with prescribed compact support.
When $\Omega$ is strictly pseudoconvex with smooth boundary, the
$\db$-closed extension
  on $\Omega$ was pioneered in the work of
Kohn and Rossi  [KoR]  using the boundary regularity of the $\db$-Neumann
operator on strongly pseudoconvex
domains.  When $\Omega \subset \Bbb C^n$ is only pseudoconvex, the
$\db$-closed extension for forms on
$M = b\Omega$ to
$\Omega$ can also be obtained if one can obtain the  boundary regularity
for the $\db$-Neumann problem or weighted
$\db$-Neumann problem (see   \cite{Sh,BSh} or   \cite{CS, Chap.9}).

In
our case,
 the connected real hypersurface $M$ divides $\Bbb CP^n$ into
two connected domains: $\Omega_+$ and $\Omega_-$. Both $\Omega_+$ and
$\Omega_-$ are pseudoconvex and  have Levi-flat
boundaries.
Using
the boundary regularity of the $\db$-Neumann
operator, a two-sided $\db$-closed extension $\tilde{\theta}_\pm$  will be
constructed  for any   form $\theta_b$
satisfying the compatibility condition (0.5) or (0.6) on the
Levi-flat hypersurface
$M$, where the $\db$-closed extension $\tilde \theta = \tilde{\theta}_\pm$
is defined on the whole space $\Bbb CP^n = \overline{ \Omega_+ \cup
\Omega_-}$.

We let   $N = N_{(p, q)}  $ denote  the
$\db$-Neumann operator on $\Omega$, i.e., the inverse operator
for the Laplace operator $\square_{(p, q)}$ acting on $(p,q)$-forms
with {\it  the $\db$-Neumann boundary condition}  on $\Omega$.
  Let
$W^{s}_{(p, q) }(\Omega)$ be the space of all $(p, q)$-forms whose
coefficients
are in the Sobolev space $W^s(\Omega) $.
To prove
Theorem 1, we  first show that $N = N_{(p, q)}|_{\Omega}$ exists and
that the operators $N$, $\db^*N$ and
 $\db N$ are bounded on $ W^s $    for some $s>0$,  when $\Omega \subset
\Bbb CP^n$ is a
pseudoconvex domain with $C^2$ boundary.

 To formulate our main  regularity result for the $\db$-Neumann operator
$N$,
we introduce a geometric invariant
for any pseudoconvex domain $\Omega \subset \Bbb CP^n$  as follows:

\proclaim{Definition 0.1} Let $\Omega$ be a  pseudoconvex  domain in a
K\"ahler manifold with $C^2$-smooth
boundary. Let $\delta(x) = \delta_{b\Omega}(x) =
d(x, b
\Omega)$ be a distance function from the boundary. We call  $t_0 =
t_0(\Omega)$ the
order of plurisubharmonicity
for
the distance function $\delta_{b\Omega  }$  if
$$ t_0(\Omega) = \sup \{ 0<\epsilon \le 1 | i\partial\bar{\partial} (-
\delta^{\epsilon   }) \ge 0 \text{ on } \Omega\}.\tag 0.7$$
\endproclaim
When $t_0(\Omega) =1$, the condition implies that there exists a
plurisubharmonic
defining function on $\Omega$. For a general smooth pseudoconvex domain  in
$\Bbb C^n$,
 such a plurisubharmonic defining function does not necessarily exist
(cf. \cite{DF2}). However,
Diederich-Fornaess \cite{DF1} showed that there exists a
$ 0<t_0(\Omega) \le 1$ for any pseudoconvex domain
$\Omega$   in
$\Bbb C^n$ with
$C^2$-smooth  boundary (for a simple proof for  pseudoconvex domains  in
$\Bbb C^n$ with $C^3$ boundary,  see Range
\cite{R}).    In
$\Bbb CP^n$,    Ohsawa-Sibony [OS] showed  that there exists $ 0<
t_0(\Omega)\le 1  $ for any
pseudoconvex domain $\Omega \subset \Bbb
CP^n$ with $C^2$-smooth  boundary using results of Takeuchi [Ta] and
[DF1].

\proclaim{Theorem 2} Let $\Omega$ be a pseudoconvex domain   with
$C^2$-smooth boundary in
$ \Bbb CP^n$ and  let $ t_0(\Omega)$ be the order of plurisubharmonicity
for the distance function $\delta_{b\Omega}$.
Then the $\db$-Neumann operator
$N_{(p, q)}$    exists on
$L^2_{(p, q) }(\Omega)$ where $0\le p,  q\le n$.   Furthermore,  $N,
\bar{\partial}N,
\bar{\partial}^*N$ and the Bergman projection $   P$ are  exact regular
  operators on
$W^{s}_{(p, q) }(\Omega)$ for $0 < s <
\frac 12 t_0 $ with respect to the  $W^s(\Omega)$-Sobolev norms.
\endproclaim

The Bergman projection $P$ is  the orthogonal projection from
$L^2_{(p, q) }(\Omega)$ to $ker(\db)$.
 In general,  $  P$ or the $\db$-Neumann operator $N$ is     not
necessarily a bounded
operator  from  $ W^s_{(p, q)}(\Omega)
 $ to $ W^s_{(p, q)}(\Omega)$ for  smooth pseudoconvex
domains $\Omega$  (cf. [Ba]).
In fact,
 Barrett [Ba] showed  that, for any given $\beta$ with $0<\beta <1$, there is a
 pseudoconvex domain $\Omega_\beta$ (the so-called Diederich-Fornaess'
worm domain \cite{DF2}) with $C^\infty$-smooth boundary
 such that the Bergman projection $  P$ and  the $\db$-Neumann operator $N$ on
$\Omega_{\beta}$ are {\it  not bounded } from $W^\beta_{(0, 1)
}(\Omega_{\beta}  )$ to
itself.    Christ [Chr]     showed that $  P $ does {\it not} map
$C^\infty(\bar{\Omega}_\beta) $ into
$C^\infty(\bar{\Omega}_\beta) $ for such a domain $\Omega_\beta$.  Therefore, one
can only expect  the boundary regularity of the
$\db$-Neumann operator
$N$  to be regular in $W^{s}_{(p, q) }(\Omega)$ with small $s$ for general
pseudoconvex domains.

When   $\Omega \subset\subset
\Bbb C^n$ with $C^\infty$-smooth boundary and $\Omega$ has a plurisubharmonic
defining function,   Boas and Straube [BS] prove
that the conclusion of  Theorem 2 holds for any
$s > 0$. This corresponds to the case when  $ t_0(\Omega) = 1$.
   For the
case
$
t_0(\Omega) < 1$ and
$\Omega
\subset
\Bbb C^n$, Kohn [Ko2] and Berndtsson-Charpentier [BC] obtain the Sobolev
regularity for
the operators $ \bar{\partial}^*N$ and the Bergman projection $  P$.
Boas-Straube [BS1] show  that
regularity for the
Bergman projection and the $\db$-Neumann operator are equivalent for
$C^\infty$-smooth pseudoconvex domains in
$\Bbb C^n$ (and also  for smooth pseudoconvex domains in complex manifolds with
strongly plurisubharmonic   functions in a
neighborhood of the boundary).  For
 pseudoconvex domains $  \Omega \subset \Bbb CP^n $,  such a
smooth  strictly plurisubharmonic function $\phi$ might not exist in a
neighborhood of the boundary
$b \Omega$. Therefore, the proof of   Theorem 2 has extra difficulties
when compared
to the case of domains
in  $\Bbb C^n$.  We remark that our proof of the regularity for the Bergman
projection and $\db^*N$ is similar to the proof
in \cite{BC}, but
the regularity results of the $\db$-Neumann operator $N$ and $\db N$ in
Theorem 2
are new, even for   domains with
$C^2$-smooth boundaries in
$\Bbb C^n$.

The    $W^s(\Omega)$-regularity result for $\db N$
     is sufficient and  crucial in the proof of
  Theorem 1. The key observation is
that  the
$\db_b$ equation (0.4), when restricted to  each complex foliation leaf of
$M$, is elliptic. Thus, any solution for (0.4)
has    H\"older regularity on each leaf from the elliptic theory.  To prove
continuity of the solution for equation
(0.4), one only needs to prove  continuity of the solution in the
transversal direction, which can be proved by using the Besov norms and a
finite difference scheme.  Details are given in
Section 5.

It should be pointed out that there exist {\it  non-smooth } Levi-flat real
hypersurfaces in
$\Bbb CP^2$.  Recall that  a  $C^1$ or Lipschitz real
 hypersurface   $M^{2n-1}$ is called  Levi-flat in a complex manifold
$\Bbb CP^{ n}$ if $\Bbb CP^n\setminus M$ consists of two pseudoconvex
domains.  For example, let $$ M=\{[z_0,z_1,z_2]\in\Bbb CP^2\mid
|z_0|=|z_1|\},$$
where $(z_0,z_1,z_2)$ are homogeneous coordinates in $\Bbb CP^2$.
Then $M$ is Levi-flat since $\Bbb CP^2\setminus M$ consists of two
pseudoconvex domains. The hypersurface $M$
is   smooth except at $(0,0,1)$, where $M$ is neither $C^1$ nor Lipschitz (in
the sense of a Lipschitz graph locally).
It is still an open question if there exist Lipschitz or $C^2$  Levi-flat
hypersurfaces in $\Bbb CP^n$.

In addition, the positive curvature condition  (or the
irreducible condition)
is needed for the possible  generalization of Theorem 1. One can inspect
the example
$S^1 \times \Bbb CP^1$, which is a $C^\infty$-smooth
Levi-flat hypersurface in $\Bbb CP^1 \times \Bbb CP^1$. Clearly,
$\Bbb CP^1 \times \Bbb CP^1$ does not  admit a K\"ahler metric with
positive holomorphic bisectional curvature, (e.g., cf. [HSW]).
Any compact K\"ahler manifold with positive  holomorphic bisectional curvature
must be biholomorphic to $\Bbb CP^n$ (e.g., cf. [Mok]). Although results in
this paper
are stated for domains in $\Bbb CP^n$, the conclusions of these results
remain to be true for corresponding domains in any compact K\"ahler
manifold with positive holomorphic bisectional curvature as well.


   The plan of this paper is as
follows: In Section 1 we
derive some   preliminary     results  on equidistant real hypersurfaces
in $\Bbb CP^n$ for the proof  of Theorems 1
and 2.   In Section 2 we discuss the
$L^2$ existence theorems of $\db$ and  the
$\db$-Neumann problem on domains  in $\Bbb CP^n$. In Section 3 we prove
that the $\db$-Neumann operator on pseudoconvex
domains with
$C^2$ boundary is regular in some Sobolev spaces as stated in Theorem 2.
Using the
existence and regularity of the
$\db$-Neumann operator,  one can study  the solution of the
$\db$-equation  with prescribed support, including the
case of forms with top degrees.  We
construct $\db$-closed extensions
from the boundary and prove the existence and regularity of equation (0.4)
in   Section 4.
 The  proof of  Theorem 1 is completed  in Section 5.

\heading   {\bf 1. Equidistant hypersurfaces and preliminaries for Theorems 1-2
 } \endheading

 In this section, we first study the geometry of  equidistant  real
hypersurfaces which will be needed to prove
Theorem 1. Afterwards we
recall four preliminary results which are indispensable  in  the proof of
Theorem  2.

   Suppose that $\Omega$ is a bounded domain with $C^2$-smooth boundary
$b\Omega$ in $\Bbb CP^n$.  Let $b\Omega_s
= \{ x \in \Bbb CP^n | d(x, b\Omega) = s\}$ and $\{b\Omega_s\}_{s <
\epsilon}$ be a
family of equidistant real hypersurfaces from $b\Omega$ in $\Bbb CP^n$.

We shall show that the connection and the curvature form  of the complex line
bundle of equidistant hypersurfaces  $\{b\Omega_s \}_{s < \epsilon}$ are
$C^{\alpha}$,  as long as $b\Omega$
is a $C^{2, \alpha}$-smooth hypersurface in $\Bbb CP^n$. We will use the
Cartan-Chern theory to prove this assertion as
follows.

Observe that, if $M = b\Omega$ is a $C^2$-smooth connected   hypersurface
in $\Bbb CP^n$, then
$M$ divides
$\Bbb CP^n$ into two connected components $\Omega_+$ and $\Omega_-$.  Let
$\rho$ be a
$C^2$ defining function for $M$ such that
$$\rho(x)=\cases &-d (z,M),\quad z\in \Omega_-,\\& d(z,M),\quad z\in
\Omega_+,\endcases\tag 1.1$$
where $d(z,M)$ denotes the distance from $z$ to $M$.
The defining function $\rho$ is called a signed distance function.
Let
$U_\epsilon(M)=\{z\in \Bbb CP^n\mid d(z,M)<\epsilon\}$ be a small tubular
neighborhood of $M$.    We    consider the real two-dimensional normal
plane bundle
$$
\Cal N_{\{\nabla\rho,J\nabla
\rho\}}=\{(c_1 \nabla\rho+c_2 J \nabla \rho)|_Q\mid Q\in U_\epsilon(M),\
c_1,c_2\in \Bbb R\}
 $$
and
the complex normal line bundle
$$
\Cal N_{ \widetilde{\nabla\rho}  }= \Cal N_{\{\nabla\rho\otimes \Bbb C\}}=\{
\lambda(\nabla\rho-\sqrt{-1}J(\nabla\rho))|_Q\mid Q\in U
_\epsilon(M),\lambda\in
\Bbb C\}. $$

The Levi-Civita connection $D$ of $\Bbb CP^n$ with the standard
Fubini-Study metric $\omega$
satisfies
$$
D _\xi(J\eta)=J(D_\xi\eta)\tag 1.2
$$
since the metric  $\omega$ is K\"ahler. We shall study the induced connections and
curvature on
$\Cal N_{\{\nabla\rho,J\nabla
\rho\}}$ and $\Cal N_{\{\nabla\rho\otimes \Bbb C\}}$. There is also a hermitian
connection $\tilde D$ on $\Bbb CT(\Bbb
CP^n)=T^{1,0}(\Bbb CP^n)\oplus T^{0,1}(\Bbb CP^n)$. The two connections $D$
and $\tilde D$ are related by
$$
\tilde D_{\xi+ {\sqrt{-1}} \eta}(Z+ {\sqrt{-1}} W)=(D_\xi Z-D_\eta W)+
{\sqrt{-1}}(D_\xi W+D_\eta Z)
$$
for $\xi,\eta,Z,W\in [T(\Bbb CP^n)]_{\Bbb R}$.

Let us consider the corresponding real orthogonal decomposition
$$
 \[T(\Bbb CP^n)\]_{\Bbb R}= \Cal N_{\{\nabla\rho,J\nabla
\rho\}}\oplus \Cal N_{\{\nabla\rho,J\nabla
\rho\}}^{\perp}.\tag 1.3
$$
For $\xi\in \[T(\Bbb CP^n)\]_{\Bbb R},$
$\xi^N $ is the projection of $\xi$ into
the normal component $\Cal N_{\{\nabla\rho,J\nabla
\rho\}}$. Similarly, we have
$$
T^{1,0}(\Bbb CP^n) = \Cal N_{\{\nabla\rho\otimes
\Bbb C\}}\oplus
\Cal N_{\{\nabla\rho\otimes\Bbb C\}}^{\perp} .\tag 1.4$$
For $\tilde\xi\in  T^{1,0}(\Bbb CP^n) ,$    ${\tilde\xi}^N $ is the
projection of  ${\tilde\xi}$ into
the complex normal component $\Cal N_{\{\nabla\rho\otimes\Bbb C\}}$.

The map $$u \mapsto  \tilde u   = \frac 1{\sqrt{2}} (u-\sqrt{-1}J u )$$ is a
linear isomorphism from $T(\Bbb CP^n)_{\Bbb R}$ to $T^{1,0}(\Bbb CP^n)$.
 We define
$$
D^N_\xi\eta=(D_\xi\eta)^N, \qquad \eta\in
\Cal N_{\{\nabla\rho,J(\nabla\rho)\}}. \tag 1.5
$$
The induced curvature tensor is defined by
$$R^N(Z,W)\eta=-D_Z^ND_W^N\eta+D_W^ND_Z^N\eta+D_{[Z,W]}^N \eta,$$
for $\eta\in N_{\{\nabla\rho,J(\nabla\rho)\}}$ and  $Z, W\in [T(\Bbb
CP^n)]_{\Bbb R} $.
Similarly, we have  the induced curvature tensor $\tilde \Theta^N$ for the
complex line bundle $ \Cal N_{
\widetilde{\nabla\rho}  }
$ on $U_\epsilon$.

 Since $\|\nabla \rho\|^2=1$,
$$d(\|\nabla\rho\|^2)=2\langle D \nabla\rho,\nabla\rho\rangle=0.$$
Thus $D_\xi^N(\nabla\rho)$ has no component in the   $\nabla\rho$ direction
and
$$
D_\xi^N(\nabla\rho)=\beta(\xi)J(\nabla\rho) \tag 1.6 $$
where  $\beta$ on $U_\epsilon$ is defined by $$\beta(\cdot
)=\text{Hess}(\rho)(\cdot
,J(\nabla\rho)).\tag 1.7
$$  The 1-form $\beta$ is the connection form for
$\Cal N_{\{\nabla\rho,J(\nabla\rho)\}}.$
Recall that from (1.2),
$D_\xi(J(\nabla\rho))=JD_\xi(\nabla\rho)$ and
$\langle
\nabla\rho,J(\nabla\rho)\rangle\equiv  0.$
Thus we have
$$D_\xi^N(J(\nabla\rho))=JD_\xi^N(\nabla\rho)=J(\beta(\xi)J(\nabla\rho))=-\beta(
\xi)\nabla\rho.\tag
1.8 $$
We shall study the connection form for $\Cal N_{\widetilde{\nabla\rho}}$
and its
curvature form $\tilde\Theta^N$.

\proclaim{Proposition 1.1} Let $M$ be a $C^m$ hypersurface in $\Bbb CP^n$
and let  $\rho$ be the signed
distance function for $M$, $m\ge 2$.  Let $\tilde e_n=\widetilde{\nabla
\rho}$ be the unit complex normal on a
tubular neighborhood $U_\epsilon(M)$ for sufficiently small $\epsilon>0$.
Then we have
 $$
\tilde D_\xi \tilde e_n=\theta(\xi)\tilde e_n \tag 1.9
$$
where
$$
\theta(\xi)=\sqrt{-1} \text{Hess}(\rho)(\xi, J(\nabla\rho)) =
\sqrt{-1} \beta(\xi).\tag
1.10
$$ Consequently,
$\sqrt{-1}\theta$ is real-valued and
the curvature form $\tilde \Theta^N$ of $\Cal N_{{ \widetilde{\nabla\rho}  }}$
is an  exact  2-form on
$U_\epsilon(M)$:
 $$\tilde\Theta^N=-d\theta. \tag 1.11 $$
Furthermore,  the connection $\theta$ and the curvature  form $\tilde\Theta^N
= - d\theta$ above are of class
 $C^{\alpha}$ if $M$
is a $C^{2, \alpha}$-smooth hypersurface in $\Bbb CP^n$.
\endproclaim
\demo{Proof}
The complex line bundle $\Cal N_{{ \widetilde{\nabla\rho}  }}$ is spanned
by  a global nowhere vanishing section
$  \widetilde {\nabla\rho}   =\frac 1{\sqrt
2}(\nabla\rho-\sqrt{-1}J(\nabla\rho))$, the complex normal vector on
$U_\epsilon(M)$.    A
direct computation shows that, by (1.6)
and (1.8) we have
$$\aligned \tilde D_\xi^N\tilde e_n
&=\frac 1{\sqrt
2}\{\tilde D_\xi^N\nabla\rho-\sqrt{-1}\tilde D_\xi^N(J(\nabla\rho))\}
=\frac 1{\sqrt
2}\{\beta(\xi)J(\nabla\rho)+\sqrt{-1}\beta(\xi)(\nabla\rho) \}
\\&=\frac {\sqrt{-1}} {\sqrt
2} \beta(\xi)\{(\nabla\rho)-\sqrt{-1}J(\nabla\rho)
\}=\sqrt{-1}\beta(\xi)\tilde e_n
=\theta(\xi)\tilde e_n,
\endaligned $$
where $\theta=\sqrt{-1}\beta$ is the connection 1-form of the complex line
bundle $N_{\widetilde{\nabla\rho}
 } $. This proves (1.9).
Using the Chern formula  \cite{Chern},  the curvature form $\tilde \Theta^N$ is
given by
$$\tilde \Theta^N=-d\theta+\theta\wedge\theta=-d\theta=-\sqrt{-1}d\beta $$
which proves (1.11).

By (1.10), $\theta$ is $C^{m-2, \alpha}$ if $M$ is $C^{m, \alpha}$.
It remains to prove that the curvature form $\tilde \Theta^N$ of $N_{{
\widetilde{\nabla\rho}  }}$ is a
 2-form with
$C^{\alpha}$ coefficients on
$U_\epsilon(M)$ if $M$ is
$C^{2, \alpha}$-smooth. This can be proved by the generalized Gauss-Codazzi
equation
(the Cartan-Chern-Gauss structure equation).
 Let $\tilde e_1,\cdots,\tilde
e_n$ be an orthonormal local basis for $T^{1,0}(U_\epsilon\cap V)$ such
that $\tilde e_n=\widetilde
{\nabla\rho}  $, where $V$ is some neighborhood near a boundary point $p\in
M$. Let $\tilde
\theta=(\tilde \theta_{k\bar l})$ be the connection 1-form  for $\Bbb CP^n$
defined
by
$$\tilde D_\xi \tilde e_j=\sum_{l=1}^n\tilde \theta_{j,\bar l}\tilde e_l. $$
We remark that though $T^{1,0}(\Bbb CP^n)$ is a holomorphic vector bundle,
our unitary frame $\{\tilde
e_1,\cdots,\tilde e_n\}$ is not necessarily holomorphic. Thus $\tilde
\theta_{k,\bar l}$ is not necessarily of
$(1,0)$ type.
 Let
 $\tilde R(X, Y) = - \tilde D_X \tilde D_Y + \tilde D_Y  \tilde D_X
+   \tilde D_{[X, Y]}$ and
$$\tilde\Theta_{k,\bar l}(Z,\overline W)=\langle \tilde R(Z,\overline
W)\tilde e_k,\bar{\tilde
e}_l\rangle,\quad Z,W\in T^{1,0}(V).\tag 1.12$$
 The matrix-valued 2-tensor
$\tilde \Theta$ is given by   Chern's formula for the  curvature of complex
vector bundles (cf. [Chern], [KN] or [Zh]),
$$\tilde \Theta=-d\tilde \theta+\tilde\theta\wedge\tilde\theta.\tag 1.13$$
In particular, we have
$$\aligned \tilde \Theta_{n,\bar n}&=\langle \tilde R(\cdot,\cdot)\tilde
e_n,\bar{\tilde e}_n\rangle
=-d\tilde\theta_{n,\bar n}+\sum_{n=1}^n\tilde \theta_{n,\bar
l}\wedge\tilde\theta_{l,\bar
n}
\\&=\tilde\Theta^N+\sum_{n=1}^{n-1}\tilde \theta_{n,\bar
l}\wedge\tilde\theta_{l,\bar
n},\endaligned\tag 1.14$$
where we have used the fact that $\tilde\theta_{n,\bar n}=\theta$ and (1.11).
It remains to calculate $\tilde \theta_{n,\bar l}$ and
$\tilde\theta_{l,\bar n}$.
For this purpose, we use the fact that
$$
\langle\tilde e_n,\bar{\tilde e}_l\rangle=\langle\tilde e_l,\bar{\tilde
e}_n\rangle=0\tag 1.15
$$
where $l=1,\cdots,n-1.$
Recall that for each $\xi\in \Bbb CT(V)$,
$$\tilde D_\xi\tilde e_n=\sum_{l=1}^m\tilde \theta_{n,\bar l}(\xi)\tilde e_l
 =\sum_{l=1}^n\langle \tilde D_\xi (\tilde e_n),\bar{\tilde e}_l\rangle
\tilde e_l.$$
Thus, we have
$$
\aligned \tilde \theta_{n,\bar l}(\xi)&=\langle \tilde D_\xi (\tilde
e_n),\bar{\tilde e}_l\rangle
=\frac 12\langle \tilde D_\xi (\nabla\rho-\sqrt{-1}J(\nabla\rho)),
e_l+\sqrt{-1}Je_l \rangle
\\&=\frac 12\langle \tilde D_\xi (\nabla\rho ),
e_l+\sqrt{-1}Je_l \rangle
\ -\frac 12\langle  \sqrt{-1}J \tilde D_\xi  (\nabla\rho) ,
e_l+\sqrt{-1}Je_l \rangle
\\&= \langle \tilde D_\xi (\nabla\rho ),
  e_l   \rangle+ {\sqrt{-1}} \langle \tilde D_\xi (\nabla\rho ),
 Je_l \rangle
=\sqrt 2\text{Hess}(\rho)(\xi,\tilde e_l).\endaligned\tag 1.16$$

Similarly, we have
$$
\aligned \tilde \theta_{l,\bar n}(\xi)&=\langle \tilde D_\xi (\tilde
e_l),\bar{\tilde e}_n \rangle
=-\langle  \tilde e_l ,\tilde D_{ \xi} (\bar{\tilde e}_n) \rangle
\\&= -\langle  \tilde D_{\bar \xi} ( {\tilde e}_n), \bar{\tilde e}_l \rangle
=-\sqrt 2\text{Hess}(\rho)(\bar\xi,\tilde e_l).\endaligned\tag 1.17$$
This implies that $$\tilde \theta_{l,\bar n}=-\bar{\tilde\theta}_{n,\bar
l},\quad\text{  on }V.\tag 1.18$$  Thus
 $$ \tilde \theta_{n,\bar l}\wedge\tilde \theta_{l,\bar n}=-\tilde
\theta_{n,\bar l}\wedge\bar{\tilde
\theta}_{n,\bar l}\tag 1.19 $$ is a skew-hermitian 2-form.

 If $\rho$ is $C^{m, \alpha}$-smooth, then
each $\tilde e_k$ is $C^{m-1, \alpha}$-smooth for $k=1,\cdots n$. Thus
$\tilde\Theta_{n,\bar n}$ is   $C^{m-1, \alpha}$-smooth.
It follows from (1.14) and (1.16)-(1.18) that
$\tilde\Theta^N$ is  a
 $C^{m-2, \alpha}$ smooth form, since $\text{Hess}(\rho)$ is $C^{m-2,
\alpha}$. The
proposition is proved. \qed
\enddemo

    The $C^{\alpha}$-regularity of $\theta$ and $\tilde\Theta^N$ will
be used in the proof of Theorem 1, see Section 5 below.
When $M = b\Omega$ is  {\it a
Levi-flat real-hypersurface},  the restriction $\Theta_b$ on $M$  of the
curvature form $\tilde \Theta^N$ above has
some additional properties.

\proclaim{Proposition 1.2}Suppose that  $M$ is a $C^2$-smooth Levi-flat
hypersurface in
$\Bbb CP^n$ with    the signed
distance function $\rho$.  Let $ \widetilde{\nabla \rho}$ be the unit
complex normal on a
tubular neighborhood $U_\epsilon(M)$ for sufficiently small $\epsilon>0$
and let $\tilde \Theta^N$ be the
curvature tensor for the complex line bundle $N_{\widetilde{\nabla \rho}}$.
 If $J$ is the complex structure of
$\Bbb CP^n$,  then $\Theta_b = \tilde{\Theta}^N|_{[T(M)]_{\Bbb R}
\cap J [T(M)]_{\Bbb R} }$ is  a  $(1, 1)$-form.
Furthermore, we have
 $$\sqrt{-1}\tilde \Theta^N(\tau,\bar\tau)  =
\sqrt{-1}\Theta_b(\tau,\bar\tau)\ge 2 \tag 1.20$$
where $\tau\in T^{1,0}(M)$ with $|\tau|=1$.
\endproclaim
\demo{Proof} If $M^{2n-1}$ is a Levi-flat
real hypersurface, then by definition we have
$$
  \text{Hess}(\rho) (\bar{\tilde{Y}},  {\tilde{X}}) =
(i \partial \db \rho)  (\bar{\tilde{Y}},  {\tilde{X}})     =  0
\tag1.21
$$
for any pair  $\{\tilde{X}, \tilde{Y}\} \in T^{(1, 0)}(M)$.

It follows from (1.16)-(1.17) and
(1.21)
that  the
connection forms
$\{ \tilde
\theta_{n,\bar l}\}$ are
$(1,0)$-forms on the Levi-flat hypersurface
$M = b\Omega$.  Using (1.14) and (1.18), we have
$$  \tilde\Theta^N= \tilde \Theta_{n,\bar n}-\sum_{n=1}^{n-1}\tilde
\theta_{n,\bar l}\wedge\tilde\theta_{l,\bar
n} =  \tilde \Theta_{n,\bar n}+\sum_{n=1}^{n-1}\tilde \theta_{n,\bar
l}\wedge\bar{\tilde\theta}_{n,\bar
l}.  \tag1.22$$
Because the Fubini-Study metric $\omega$ is a K\"ahler metric, its curvature
form $\tilde \Theta$ is a
matrix valued (1, 1)-form with respect to any unitary frame (cf. [Chern],
[Zh]). Hence,
$ \tilde \Theta_{n,\bar n}$ is a  (1, 1)-form.
Since   $\{ \tilde
\theta_{n,\bar l}\}$ are
$(1,0)$-forms on the Levi-flat hypersurface
$M = b\Omega$, we see that
 $\sqrt{-1}\sum_{n=1}^{n-1}\tilde \theta_{n,\bar
l}\wedge\bar{\tilde\theta}_{n,\bar
l}$ is a nonnegative $(1,1)$-form.  This fact together with (1.22) implies
that
$\Theta_b $ is a $(1,1)$-form.  Furthermore, it follows that, for $\tau \in
T^{(1,0)}(M)$ with $|\tau| = 1$,
$$\sqrt{-1} \tilde\Theta^N(\tau,\bar\tau)\ge \sqrt{-1} \tilde
\Theta_{n,\bar n}(\tau,\bar\tau)\ge 2,$$
where we used the
fact that the Fubini-Study metric of $\Bbb CP^n$ has positive holomorphic
bisectional
curvature $\ge 2$, (e.g., cf. [KN], [Zh]).
This completes the proof. \qed
\enddemo

   Let us now  recall some
preliminary results, which play important roles in the proof of Theorem 2.
First of all,
we need to recall the so-called $\db$-Neumann boundary condition. The
$L^2$-theory
of the $\db$-equation $\db u = v$ and the Laplace equation
$\square u = f$ are  related to the formal adjoint $\vartheta$
of $\db$ on the space $L^2_{(p, q)}(\Omega)$. Let   $*$ be  the
real Hodge star operator of the Riemannian metric and let  $*$ be
extended $\Bbb C$-linearly on the space of $(p, q)$-forms. Then
$$
\vartheta u = - \bar{*}[ \db (\bar{*} u)] = - * \partial * u. \tag1.23
$$
For $u,  v\in L^2(\Omega)$, we use
$$(u,\bar v)_{L^2(\Omega)}=\int_\Omega\langle u,\bar v\rangle$$
to denote the inner product.
For $q >0$, using integration by parts on the domain $\Omega$ with
$C^2$-smooth boundary
$b\Omega$, we have that $u \in \text{Dom}(\db^*|_\Omega)\cap C^1_{(p,
q)}(\overline\Omega)$ if and only if
$$
 (\vartheta u, \bar v )_{L^2(\Omega)} =  ( \db^* u, \bar v)_{L^2(\Omega)}
=  ( u, \overline{\db v} )_{L^2(\Omega)} \tag1.24
$$
holds for all $v \in L^2_{(p, q)}(\Omega)$. The condition (1.24)
is equivalent to the so-called $\db$-Neumann boundary condition
$$
   u \llcorner_{(\db \rho)_\#} = 0 \text{ on } b\Omega, \tag1.25.1
$$
where $(\db \rho)_\# = \frac 1{\sqrt{2}} (\nabla \rho + \sqrt{-1}J
\nabla \rho ) $ is the dual of $\db \rho$. In general, $u  \in
\text{Dom}(\db^*|_\Omega)$ implies that
(1.25.1) holds as currents.

   Throughout this paper, we let
$$\text{Dom}(\db^*) =\text{Dom}(\db^*|_\Omega)
= \{ u \in L^2_{(p, q)}(\Omega) |   u\llcorner_{(\db \rho)_\#} = 0 \text{
on } b\Omega  \} \tag1.25.2$$
 and
$$ \text{Dom}(\square|_\Omega) =
\{ u \in L^2_{(p, q)}(\Omega) | u\llcorner_{(\db \rho)_\#} = 0
\text{ and } (\db u)\llcorner_{(\db \rho)_\#} = 0 \text{ on } b\Omega
\}. \tag1.25.3$$

    There always exists the {\it weak} Hodge-Kodaira decomposition for any
domain $\Omega$:
$$
L^2_{(p, q)}(\Omega) =\text{ker}(\square) \oplus_{L^2}
\overline{\text{Range}(\square|_\Omega)}, \tag1.26
$$
where $\overline{\text{Range}(\square|_\Omega)}$ denotes the closure of
${\text{Range}(\square|_\Omega)}$ in $L^2(\Omega)$.
A necessary condition for the existence of the $\db$-Neumann
operator  on $L^2_{(p, q)}(\Omega)$  is the condition that the range of
$\square $ must be closed, i.e.,
$\overline{\text{Range}(\square)} = \text{Range}(\square)$. To prove
the first part of Theorem 2, by the definition of $\square = \db \db^* +
\db^* \db$,  it is sufficient to
 show that both $\text{Range}(\db)$ and $\text{Range}(\db^*)$ are closed
subspaces of $L^2_{(p, q)}(\Omega)$.
For this purpose (and for the convenience of the reader), we recall
the following  elementary  but very useful fact in functional analysis
(cf.  [H\"o2-3] or [CS, p60]):

\proclaim{Lemma 1.3}  Let $\Cal F: H_1
\to H_2$ be a linear, closed, densely defined operator between two Hilbert
spaces.
The following conditions on $\Cal F$ are equivalent:
\medskip
\noindent
(1) $\text{Range}(\Cal F)$ is closed;
\medskip
\noindent
(2) There is a constant $C$ such that
$$
     \| u \|_1 \le C \| \Cal F u \|_2 \text{ for all } u \in \text{Dom}(\Cal F)
\cap \overline{\text{Range}(\Cal F^*)}; \tag1.27.1
$$
\medskip
\noindent
(3) $\text{Range}(\Cal F^*)$ is closed;
\medskip
\noindent
(4) There is a constant $C$ such that
$$
     \| v \|_2 \le C \| \Cal F^* v\|_2 \text{ for all } v \in
\text{Dom}(\Cal F^*)
\cap \overline{\text{Range}(\Cal F)}; \tag1.27.2
$$
The best constants in (1.27.1)-(1.27.2) are the same.
\endproclaim

   Lemma 1.3 will be used in the proofs of Theorem 2.1, Theorem 2.6
and Theorems 3.4-3.5 below, which are parts of Theorem 2. By Lemma 1.3,
it suffices  to derive {\it a priori} estimates for the
$\db$ equation and the $\square$ equation.
We recall the curved version of the Morrey-Kohn-H\"ormander formula,
which is a Bochner type formula with a weight function $e^{ -\phi}$.
The weighted function $e^{-\phi}$ induces  a perturbed adjoint operator:
$$
   \vartheta_\phi u = \db^*_\phi u = e^\phi \vartheta (e^{-\phi u}) =
- e^\phi * \partial (e^{-\phi}*u). \tag 1.28
$$
Of course, one can formulate the  weighted Laplace operator
$\square_\phi  = (\db \db^*_\phi + \db^*_\phi \db)  $ and the
corresponding weighted $\db$-Neumann operators $N_{\phi}$.

     Our strategy is to use the weighted $\db$-Neumann operators $ N_{\phi}$ to
 estimate the original $\db$-Neumann operator $N$, see Sections 2-3 below.
In order to estimate the operator $\db^*_\phi N_\phi$, we
let $L_1,\cdots,L_n$ be a local orthonormal frame for $T^{(1,0)}(U)$, where
$U$ is
a local neighborhood near some point in $\overline \Omega$. For any $u\in
C_{(p,q)}(\overline\Omega)$ and $\phi\in  C^2(\Omega)$, we define an operator
$$ \langle( i\partial\db \phi) u, \overline u\rangle=\sum_{j,k=1}^n
(i\partial\db
 \phi) (L_j,\overline L_k)   \langle u\llcorner \overline L_j, \overline
u\llcorner
L_k\rangle \text{   on } \Omega. \tag1.29
$$
Notice that $\langle ( i\partial\db \phi)u, \overline u\rangle$ is
independent of the
choice of a local unitary basis.

Similarly,
for $Q\in b\Omega$, we require that $L_n = \sqrt{2} (\partial \rho)_\#$
and define
$$\langle  (i\partial\db \rho) u, \overline u\rangle=\sum_{j,k=1}^{n-1}
(i\partial\db
 \rho) (L_j,\overline L_k)   \langle u\llcorner \overline L_j, \overline
u\llcorner
L_k\rangle \text{    on } b\Omega. \tag1.30
$$

In addition, since $\Bbb CP^n$ has non-zero curvature,  we set
$$ \langle \Theta u, \overline u\rangle=\sum_{j,k}^n
\langle \overline
 w^j \wedge\{ [R(L_j , \overline L_k) u]\llcorner \overline L_k\},  \overline
 u\rangle   \text{   on } \Omega,
 \tag1.31
$$
where $\{w^j \}$ is the dual  of the frame $\{ L_k \}$ and $R_{X, Y}  = -
D_X D_Y u + D_Y D_X + D_{[X, Y]}$ is the curvature
operator on $(p, q)$-forms, see [Wu, Chapter 2]. Finally, we set
$$
(u, \overline v)_\phi = \int_\Omega \langle u, \overline v \rangle
e^{-\phi}. \tag1.32
$$

\proclaim{Proposition 1.4} (Bochner-H\"ormander-Kohn-Morrey formula) Let
$\Omega$ be a compact domain $\Omega$
with $C^2$-smooth boundary $b\Omega$ and $\rho(x) = - d(x, b\Omega)$.
Then, for any (p,q)-form $u\in
\text{Dom}(\db)\cap\text{Dom}(\db^*)$
 with $q\ge1$ on $\Omega$, we have
$$
\|\db u\|_\phi^2+\|\db^*_\phi
u\|_\phi^2=\|\overline\nabla
u\|_\phi^2 +(\Theta u, \overline u)_\phi
+  ((i\partial\db \phi) u, \overline u )_\phi+\int_{b\Omega} \langle
 (i\partial\db \rho) u, \overline u\rangle e^{-\phi}, \tag1.33
$$
where $ \|\overline\nabla
u\|_\phi^2 = \int_\Omega \sum^n_{j = 1} |D_{\bar L_j} u|^2 e^{-\phi}$
and $\{ L_1, ..., L_n\}$ is a local unitary frame of $T^{(1,0)}(\Omega)$.
\endproclaim
\demo{Proof} This formula is  known (cf. [AV], [H\"o2]
 [Dem1-2], [Siu1],    [Wu]) for some special cases, although it has not
been stated in the literature in the
form (1.33).   If $u$ has
compact support in the interior of $\Omega$, the formula (1.33) was proved
in  [AV], Chapter 8 of
[Dem2] and (2.12) of [Wu].

It remains to discuss the boundary term of (1.33). For the case
$\phi = 0$, the stated formula with the boundary term was proved in [Siu1].
To compute the boundary term,
one notice that
although the boundary terms in the proof
of (1.33) involve the weight function, they have nothing to do with the
Riemannian curvature $R$ nor
$\Theta$.  The boundary term had been computed in H\"ormander
\cite{H\"o2, Chapter 3}  by combining
   the Morrey-Kohn technique
  on the boundary   with non-trivial weight function. If one combines the
results of \cite{H\"o2} with the interior formulae
discussed above,
 one can prove that
  (1.33) holds for the general case with a weight function $e^{-\phi}$
and the  curvature term.
   \qed
\enddemo

In later application,  we will choose $\phi = - t \log |\rho|$ in (1.33).
Therefore, we need to discuss the term $ i \partial \db \phi$ with
$ \phi = - t \log |\rho|$ and the curvature term
$\langle \Theta u, \bar u \rangle$ in the right hand side of (1.33).

\proclaim{Proposition 1.5}  Let $u$ be a $(p, q)$-form of $\Omega \subset
\Bbb CP^n$ with $q \ge 1$. Suppose that $\Theta$ is the curvature
term defined in (1.31) with respect to the Fubini-Study metric. Then
 we have
  $$  \aligned \cases
   \langle \Theta u, \overline u\rangle &= 0,\quad \text{for any $(n,
q)$-form u;}  \\      \langle
\Theta u,
\overline u\rangle &\ge 0, \quad\text{ when $p \ge 1$ and $u$ is a (p, q)-form}
 ;\\    \langle \Theta u, \overline u\rangle  &= q(2n+1)
|u|^2, \quad\text{ when $u$ is a (0, q)-form}.\endcases
\endaligned\tag 1.34$$

\endproclaim
\demo{Proof} The assertion for  $(0, q)$-forms and
$(n, q)$-forms was computed in
[Wu] and [Siu1]. For the curvature operator $\Theta$ acting on (1, 1)-forms
$u$,
[Siu1] showed that $\Theta$ is non-negative definite, but {\it not}
positive definite,
(see [Pe] as well). Henkin and Iordan [HI] extended this observation for all
$(p, q)$-forms with $p \ge 1$. In fact,
Lemma 3.3 of [HI] and its proof showed that $\Theta$ acting on
$L^2_{(p, q)}(\Omega) $ is a non-negative operator. Furthermore, $ \Theta$ is
{\it positive definite}
on $L^2_{(p, q)}(\Omega) $ if and only if $p=0$ and $q \ge 1$.
  \qed\enddemo

\proclaim{Proposition  1.6}
 If $ \Omega$ is a   pseudoconvex domain with
$C^2$-smooth boundary $b\Omega$ in $\Bbb CP^n$ and $\rho$ is the signed
distance function,   then
$$
i\partial\db (-\log|\rho|)(X, \overline X)
\ge \frac 12 |X|^2 \tag 1.35
$$ for any $X \in T^{(1,0)}(\Omega). $
\endproclaim
Proposition 1.6 was proved in Takeuchi
(e.g., cf. [Su], [Ta]). A new   proof of Proposition 1.6 using comparison
theorems in Riemannian geometry
will appear  in another paper in [CaS].

\heading {\bf 2. $L^2$ theory for $\db$ on pseudoconvex
domains in $\Bbb
CP^n$}\endheading
Let $\Omega$ be a bounded pseudoconvex  domain in $\Bbb CP^n$. We fix our
Fubini-Study metric $\omega$ for
$\Bbb CP^n$.
In this section, we will use preliminary results of Section 1 to study
the $\db$-Neumann operators $N$,  $\db^* N$ and $\db N$ acting on $L^2_{(p,
q)}(\Omega)$.
We first prove the  $L^2$  existence theorem of
the $\db$-Neumann operator for the easier case when $p=0$.

\proclaim{\bf Theorem 2.1}  Let $\Omega$ be a pseudoconvex domain
 in $\Bbb CP^n$, $n\ge 2$.  Then the $\db$-Neumann operator $N_{(0,q)}$
exists and is bounded on
$L^2_{(0,q)}(\Omega)$,
  where $0\le q\le n $. For any $f\in L^2_{(0,q)}(\Omega )$,
  $$\aligned f&=
\bar\partial
  \bar\partial^*N_{(0,q)}f \oplus\bar\partial^*\bar\partial
N_{(0,q)}f,\quad 1\le q\le n-1.\\
   f&=
 \bar\partial^*\bar\partial
N_{(0,0)}f\oplus Pf,\quad q=0,\endaligned\tag 2.1$$
where $P:L^2(\Omega)\to L^2(\Omega)\cap\text{Ker}(\db)$ is the Bergman
projection
with
$$P=I-\db^*N_{(0,1)}\db \tag 2.2$$
and
$$N_{(0,0)}=\vartheta N_{(0,1)}^2\db.\tag 2.3$$
Moreover, $N$, $\db^* N$ and $\db N$ are bounded operators with respect to
the $L^2$-norms.
\endproclaim
\demo{Proof}  We first assume that  $\Omega$ is a pseudoconvex domain
 in $\Bbb CP^n$  with   $C^2$
boundary $M$. Let $1\le q\le n $ and let
$\rho$ be a $C^2$ defining function for $\Omega$ normalized by $|d\rho|=1$
on $M$.  Choosing $\phi = 0$ in Proposition 1.4 , we have   for any
$(0,q)$-form $u\in
\text{Dom}(\db)\cap\text{Dom}(\db^*)$,
$$
\|\db u\|^2+\|\db^*
u\|^2=\|\overline\nabla
u\|^2+\int_M\langle (i\partial\db \rho)u, \bar  u\rangle +(\Theta u,\bar   u)
\ge \|\overline\nabla
u\|^2 + (\Theta u,\bar  u ),\tag 2.4.1
 $$ where
$\Theta$ is the Ricci form and $\overline\nabla $ is the holomorphic
gradient.   Since the Ricci curvature of $\Bbb CP^n$ with the Fubini-Study
metric is equal to
2n + 1, (cf. [KN] [Pe] or [Zh]), by (2.4.1) we have
$$
Q(u, \bar  u) = \|\db u\|^2+\|\db^* u\|^2\ge q(2n + 1)\|u\|^2,\quad u\in
\text{Dom}(\db) \cap \text{Dom}(\db^*). \tag2.4.2  $$
Consequently, we have
$$
(\square u,  \bar u) = Q(u,\bar  u ) \ge q(2n + 1)\|u\|^2,\quad u\in
\text{Dom}(\square_{(0, q)}|_\Omega),\tag2.4.3
 $$
for $q \ge 1$.

It follows from (2.4.3)  that $\text{Ker}(\square_{(0, q)}) = 0$  and $\|
N_{(0, q)}\|_{L^2} \le \frac{1}{q(2n+1)}$
 for $q \ge 1$.
Moreover, by (2.4.2)-(2.4.3) and Lemma 1.3, one can show that the three
operators $\db, \db^*$ and $\square$
have closed ranges in $L^2(\Omega)$. Notice that if $\db f = 0$ and if $u =
\db^* N f$, then $u$ is the
solution of $\db u = f$ with the smallest $L^2$-norm. Using this fact, by
(2.4.2) and Lemma 1.3, one can further show that
$\db^*N$ is a bounded operator. Since $\db N = (N\db^*)^* = (\db^* N)^*$,
using Lemma 1.3 again, we see that
the operator $\db N$ is bounded as well. This proves Theorem 2.1 for the
case of $q \ge 1$.

The existence of $N_{(0, 0)}$  also follows and one can prove (2.2) and
(2.3) exactly as for domains
in $\Bbb C^n$.  The general case for nonsmooth domains follows
from exhausting a pseudoconvex domain  by $C^2$ pseudoconvex domains and
using (2.4)  (See proofs of  Theorem
4.4.1 and Theorem 4.4.3 in   [CS] for details). \qed
\enddemo

When   $0<p\le n$, the proof of the $L^2$ existence of the $\db$-Neumann
operator is more involved.  we will use Lemma 1.3 and Propositions 1.4-1.6
with non-trivial  weight function $e^{-\phi} = |\rho|^t$.
 Here is a preliminary result on the complex Hessian of the weight function.

\proclaim{\bf Lemma 2.2} Let $\Omega\subset\subset\Bbb CP^n$
be a pseudoconvex domain with  $C^2$-smooth boundary
$b\Omega$ and let $\delta(x) = d(x, b\Omega)$ be the distance function to
$b\Omega$. Let  $t_0 = t_0(\Omega)$ be the
order of plurisubharmonicity
for
the distance function  $\delta$ defined by (0.7). Then, for any $0<t <  t_0$,
there exist    $C_t>0$ such that
$$
i\partial\db(- \delta^t)\ge C_t   \delta^t \(\omega+
\frac{i\partial\delta\wedge\db\delta }{\delta^2} \),
\tag 2.5$$
where $\omega$ is the K\"ahler form of the Fubini-Study metric on $\Bbb CP^n$.
\endproclaim
 \demo{Proof}
The existence of such $t_0$ for   pseudoconvex domains in
$\Bbb C^n$ with $C^2$-smooth  boundary   in $\Bbb CP^n$ is proved in
Ohsawa-Sibony
[OS]. Thus, by the inequality $ i\partial\db(- \delta^{t_0})  \ge 0 $
with
$0<t_0 \le 1$, we obtain  that
$$i\frac{\partial\db(-\delta)}\delta+(1-t_0)\frac
{i\partial\delta\wedge\db\delta}{\delta^2}\ge 0\tag
2.6$$
Using Proposition  1.6,  we have
$$
i\partial\db(-\log\delta)=i\frac{\partial\db(-\delta)}\delta+ \frac
{i\partial\delta\wedge\db\delta}{\delta^2} \ge  C\omega\tag 2.7
$$ where $C = \frac 12$.
Multiplying (2.6) by $(1-\epsilon)$ and multiplying (2.7) by $\epsilon$,
 and adding the two inequalities together (i.e., using the sum
$(1-\epsilon)(2.6) + \epsilon (2.7)$), we conclude  that, for any $0\le
\epsilon\le 1$,
the inequality
$$
i\partial\db(-\log\delta)= i\frac{\partial\db(-\delta)}\delta+ \frac
{i\partial\delta\wedge\db\delta}{\delta^2}\ge C\epsilon
\omega+(1-\epsilon)t_0\frac
{i\partial\delta\wedge\db\delta}{\delta^2}\tag 2.8
$$
holds.
Hence,  for any $0<t<t_0$, we choose $\epsilon_t$ such that
$(1-\epsilon_t )t_0>t$. Then
$$
\aligned
i \partial\db(- \delta^t)
&= it\delta^t \(\frac {\partial\db(-\delta)}\delta+(1-t)
\frac{ \partial\delta\wedge\db\delta }{\delta^2} \)
\\& \ge C_t  t
\delta^t
\(\omega+
\frac{i\partial\delta\wedge\db\delta }{\delta^2}\)
\endaligned
$$
where $C_t =\min( \frac 12 \epsilon_t ,(1-\epsilon_t )t_0-t)$.
The lemma is proved. \qed
\enddemo

\medskip

  We  will  use   H\"ormander's weight function method to obtain the
$L^2$ existence for the $\db$-equation with weights first.
Let $t $ be any real number and $\phi\in C^2(\Omega)$.  Let $L^2(\delta^t)$
denote the  $L^2$ space  with respect to
the weight function
$e^{-t\phi}=\delta^t$ and
$$\|f \|_{(t)}^2=\int_\Omega |f|^2e^{-t\phi}=\int_\Omega \delta^t|f|^2 .
$$
We
use $\db^*_{ t }$ to denote the adjoint of
$\db$ with respect to the weighted space.

\medskip

   We  use the norm  $\|. \|_{(t)}^2$ since  it is equivalent to
the Sobolev norm on a sub-space of $W^{-\frac t2}(\Omega)$, (see [CS p348]
and [MS1-2]).
In fact, we have
$$
\|\tilde u \|_{(t)}^2 \simeq \| \tilde u \|_{W^{-\frac t2}(\Omega ) } \tag2.9
$$
 when $\tilde u$  satisfies an elliptic equation such as $\db\oplus\vartheta$.

In order to use Proposition 1.4 to derive the desired estimates for the
operator $\db^* N$,
we need to introduce the following two asymmetric weighted norms. These new
norms will be
used to obtain    more refined   estimates than those obtained in [H\"o2].

For any $(p,q)$-form $f$ on $\Omega$, we decompose $f$ into complex normal
and tangential parts by setting
$$\cases  f^\nu&=(f\llcorner_{(\db\delta)_\#})\wedge \db\delta \\ f^\tau
&=f-f^\nu.\endcases$$
The above decomposition is well-defined for any $(p,q)$-form $f$ supported
in \newline
$[\Omega\setminus \text{Cut}_\Omega(b\Omega)]$,
where $ \text{Cut}_\Omega(b\Omega)$ is the cut-loci of $b\Omega$ in
$\Omega$,
see [Cha], [CE] or [Pe]. It is well-known that
 the cut-loci $\text{Cut}_\Omega(b\Omega)$ of $b\Omega$ has real
dimension $\le (2n-1)$ and hence has zero measure in
$\Omega$,
 (e.g., [Cha], [CE, p90ff], [Pe]). In summary, the above decomposition
exists {\it almost everywhere in
$\Omega$}.

We define the asymmetric weighted norm
$$|f|^2_A=|  f^\tau|^2+ \frac{|  f^\nu|^2}{|\delta|^2}.\tag 2.10.1$$

We also define the dual norm
$$|f|^2_{A'}=|  f^\tau|^2+  {|  f^\nu|^2}{|\delta|^2}.\tag 2.10.2$$
For any $t>0$, let $
L^{2}_A(\delta^t)$ and $  L^{ 2}_{A'}(\delta^t)$ denote the
weighted $L^2$ spaces  on $(p,q)$-forms defined by the norm
$$|||u|||_{(t)}^2=\int_{\Omega}\delta^t |f|_A^2=\int_{\Omega}\delta^t( |
f^\tau|^2+ \frac{|  f^\nu|^2}{|\delta|^2}) \tag 2.10.3 $$
and
$$|||u|||_{ (t)'}^2=\int_{\Omega}\delta^t |f|^2_{A'}=\int_{\Omega}\delta^t
| ( f^\tau|^2+  {|
f^\nu|^2}{|\delta|^2}). \tag2.10.4$$

We may assume that $\text{Diam}(\Bbb CP^n) \le 1$ up to a factor
$\frac{\pi}{2}$. It is  obvious that
$$|||u|||_{ (t)'}^2\le \p u\p_{ (t)}^2\le |||u|||_{ (t) }^2,$$
if $\text{Diam}(\Omega) \le 1$.

   The following is a preliminary estimate for the operator $\db^* N$,
because $u = \db^* N f$ is the solution
of $\db u = f$ with the least $L^2$-norm.

\proclaim{\bf Proposition 2.3} Let $\Omega\subset\subset\Bbb CP^n$
be a pseudoconvex domain with  $C^2$-smooth boundary
$b\Omega$.  Let  $t_0 = t_0(\Omega)$ be the
order of plurisubharmonicity
for
the distance function  $\delta$ defined by (0.7). For any  $ 0 < t < t_0$
and any $(p,q)$-form
  $f\in L^{2}_{A'}(\delta^t )
$, where
$0\le p\le n$ and
$1\le q\le n$, such that
$\db f=0$ in
$\Omega$,
   there exists $u\in L^2_{(p,q-1)}(\delta^t)$ satisfying $\db u=f$ and
$$\|u\|_{(t)}^2 \le  \frac 1{C_0t}|||f|||^2_{(t)'}. \tag 2.11  $$
\endproclaim
\demo{Proof} By  Proposition 1.6, we have that $\phi=-\log \delta$ is
strictly plurisubharmonic
and $i\partial\db\phi\ge  \frac 12 \omega$, where $\omega$ is the K\"ahler
form of $\Bbb CP^n$ with the Fubini-Study metric. It is known that
$\text{Dom}(\db^*_{t }) = \text{Dom}(\db^*)$,
(e.g., [CS, Chapter 4]).
     Using  H\"ormander's weighted estimates (cf. Proposition 1.4) and
Proposition 1.5,
 we have  the
following   formula: for any $(p,q)$-form $g\in
\text{Dom}(\db)\cap\text{Dom}(\db^*_{t })$,
$$
\|\db g\|^2_{(t)}+\|\db^*_{t
} g\|^2_{(t)}\ge t((i\partial
\db) \phi g,\bar g)_{(t)}. \tag 2.12.1
$$
From (2.8), the positive (1,1)-form
$i\partial
\db \phi$  induces a  pointwise norm on $(p,q)$ forms
$$
|g|^2_{i\partial\db\phi}= \langle i(\partial\db\phi)  g, \bar g\rangle\ge C_0
|g |_A^2 ,  \tag2.12.2
$$
where $C_0=\min{(\frac 12,t_0)}$ with  $t_0 = t_0(\Omega)$ as in (2.5).
 Thus, we have
$$\|\db g\|^2_{(t)}+\|\db^*_{t } g\|^2_{(t)}\ge C_0t|||g|||_{(t)}. \tag 2.13$$

Let $(i\partial\db\phi)'$ denote the dual norm for $(p,q)$-forms induced by
$i\partial\db\phi$.
Using (2.13) and the same argument as in  Demailly \cite{Dem1,2}, we conclude
that
for any $f  \in L^{2}_{A'}(\delta^t )$, there exist $u\in L^2(\delta^t)$
satisfying $\db u=f$ and
$$
\int_\Omega|u|^2\delta^t\le  \int_\Omega |\db u|_{(i\partial\db
\phi)'}^2\delta^t \le \frac 1{ C_0t}
\int_\Omega |\db u|_{A'}^2\delta^t.\tag 2.14
$$
This completes the proof. \qed
\enddemo

   Notice that the terms in (2.14) have an extra weight factor $\delta^t$.
We need to remove
this factor $\delta^t$   in order to obtain the desired $L^2$-estimate for
all $p \ge 0$. The next two propositions have already been obtained in
Berndtsson-Charpentier
\cite{BC} (see also \cite{HI}). We include the
proof here for the sake of completeness.

\proclaim{\bf Proposition 2.4} Let $\Omega\subset\subset\Bbb CP^n$
be a pseudoconvex domain with  $C^2$ boundary
$b\Omega$.  For   any
$f\in L^2_{(p,q)}(\Omega)
$, where
$0\le p\le n$ and
$1\le q\le n$, such that
$\db f=0$ in
$\Omega$, there exists $u\in L^2_{(p,q-1)}(\Omega)$ satisfying $\db u=f$
with $\int_\Omega|u|^2 \le \hat{C} \int_\Omega|f|^2$.
\endproclaim

\demo{Proof}
By  Proposition 2.3,   for any $t>0$ with $0 < t < t_0/4$,  there exists $u\in
L^2_{(p,q-1)}(\delta^t)$ satisfying $\db u=f$, such that
$u$ is perpendicular to Ker$(\db)$ in
$L^2(\delta^t)$
and
 $u$ satisfies inequality (2.11).

Consider $v=u\delta^{-t}$. Then $v\in L^2(\delta^{2t})$ and $v\perp
\text{Ker}(\db)$ in $L^2(\delta^{2t})$.
It follows from (2.14) that the following holds:
$$ \int_\Omega|u|^2 =   \int_\Omega|v|^2\delta^{2 t} \le \frac
1{2C_0t}\int_\Omega|\db
v|^2_{A^{'}}\delta^{2t}.
\tag 2.15  $$

   By the definitions of the norm $|  .|_{ A'}$ and $\phi = - \log \delta$
(cf. (2.10.2)), we have
$$\aligned  |\db v|^2_{A^{'}}e^{-2t\phi}&\le 2|\db
u|^2_{A^{'}}+2t^2|\db  \phi |^2_{A^{'}}|u|^2
\\&  \le   C( |\db
u|^2 +2t^2|u|^2).\endaligned\tag 2.16  $$
Choosing  $t$ sufficiently small and substituting (2.16) into (2.15),
  one obtains
  $$\int_\Omega|u|^2 \le \hat{C}\int_\Omega|\db u|^2, \tag2.17 $$
where $ \hat{C} = \frac{C}{2 C_0 t(1 - 2 C t^2   )}$.
This proves the theorem. \qed
\enddemo
By (2.9), the Sobolev norm $\| u \|_{W^{\frac t2}_{(p, q)}(\Omega)} $ of
the solution u is related to
the weighted norm $\| u \|_{(-t)}$.

\proclaim{\bf Proposition 2.5} Let $\Omega\subset\subset\Bbb CP^n$
be a pseudoconvex domain with  $C^2$ boundary
$b\Omega$. Let  $t_0 = t_0(\Omega)$ be the
order of plurisubharmonicity
for
the distance function  $\delta$ defined by (0.7).  For   any
$0\le t< t_0$,
$f\in L^2_{A'}(\delta^{-t})
$, where
$0\le p\le n$ and
$1\le q\le n$, such that
$\db f=0$ in
$\Omega$, then  $u=\db^*N f$ is in $ L^2_{(p,q-1)}(\delta^{-t})$
with $\db u=f$ and
$$ \|u\|_{(-t)}^2 \le \tilde  C_t|||f|||^2_{(-t)'}, $$
where the weighted norms are given by (2.10.1)-(2.10.4).
\endproclaim
\demo{Proof}
Let $v=u\delta^{-t}$. We have that $\db^*_t v=0$ if $q>1$ and $v\perp
\text{Ker}(\db)$ in $L^2(\delta^{t})$
for $q=0$. This implies that  $v=\db^*_tN_t\db v$.
By (2.8) (with $\epsilon=0$),  we have that
$$|\frac {\db\delta}\delta\wedge
u|^2_{(i\partial\db
\phi)'}\le
\frac 1{t_0}|\frac {\db\delta}\delta\wedge u|^2_{A'}\le \frac
1{t_0}|u|^2.\tag 2.18
$$
Thus for any $\epsilon>0$, by Proposition 2.3 (the proof of (2.14)),
we have
$$
\aligned
\int_\Omega |u|^2\delta^{-t}
&=\int_\Omega |v|^2\delta^t \le
\frac {1}{C_0t}\int_\Omega |\db v|^2_{(i\partial\db\phi)'}\delta^t
\\&\le (1+\frac 1\epsilon )\frac 1{C_0t}\int_\Omega
|\db u|_{(i\partial\db\phi)'}^2\delta^{ -t}+ (1+\epsilon)\frac{t
}{t_0}\int_\Omega |\frac
{\db\delta}\delta\wedge u|^2_{ A '}\delta^{- t}
\\&\le (1+\frac 1\epsilon )\frac 1{C_0t}\int_\Omega
|\db u|_{A'}^2\delta^{ -t}+ (1+\epsilon)\frac{t }{t_0}\int_\Omega |u|^2
\delta^{-
t}.
\endaligned\tag 2.19$$
Choosing $\epsilon$ sufficiently small such that $(1+\epsilon)\frac{t
}{t_0}=r<1$ in
 (2.19), we conclude that

$$\int_\Omega |u|^2\delta^{-t}\le \frac{ (1+\frac 1\epsilon )\frac
1{C_0t}}{1-r}  \int_\Omega
|\db u|_{A'}^2\delta^{ -t}\le \tilde{ C_t}|||f|||^2_{(-t)'}.\tag 2.20$$
Proposition 2.5 is proved. \qed
\enddemo

 We arrive at the main result of this section.
 \proclaim{\bf Theorem 2.6} Let $\Omega\subset\subset\Bbb CP^n$
be a pseudoconvex domain with  $C^2$-smooth  boundary
$b\Omega$. Then $\square_{(p,q)}$ has closed range and  the $\db$-Neumann
operator
$N_{(p,q)} : L^2_{(p,q)}(\Omega)\to L^2_{(p,q)}(\Omega)$
exists for every $ p,\ q $ such that  $0\le p\le
n, 0\le q\le n $. Moreover,
 for any $f\in L^2_{(p,q)}(\Omega )$, we have
  $$\aligned f&=
\bar\partial
  \bar\partial^*N_{(p,q)}f \oplus\bar\partial^*\bar\partial
N_{(0,q)}f,\quad 1\le q\le n-1.\\
   f&=
 \bar\partial^*\bar\partial
N_{(p,0)}f\oplus Pf,\quad q=0,\endaligned$$
where $P$ is the   projection  from $L^2_{(p,0)}(\Omega)$ onto $
L^2_{(p,0)}(\Omega)\cap\text{Ker}(\db)$ and
$$N_{(p,0)}=\db^*N_{(p,1)}^2\db.$$
In addition, $\db N$, $\db^* N$ and $N$ are bounded linear operators on
$L^2_{(p,q)}(\Omega )$.
\endproclaim
\demo{Proof}
The $L^2$-existence theorem for the $\db$-Neumann operator $N$ on $\Omega$
follows from the $L^2$-existence of the solution $u$ for the
$\db$-equation, which
is  proved
in  Proposition 2.4. The proof of Theorem 2.6 follows from Lemma 1.3
using Proposition 2.4 and (2.17)      (see e.g., proofs of
Theorem 4.4.1 and Theorem 4.4.3 in
\cite{CS} or [H\"o3].) \qed
\enddemo
We remark that when $p=0$, Theorem 2.6 holds for pseudoconvex domains not
necessarily with $C^2$ boundary from
Theorem 2.1.

\heading{\bf 3. Sobolev estimates for the $\db$-Neumann operator on
pseudoconvex domains in $\Bbb
CP^n$}\endheading

In this section  we  prove that the $\db$-Neumann operator $N$ is a {\it
bounded } linear operator on
Sobolev spaces $W^t_{(p, q)} (\Omega)$ for small $t>0$.
In order to derive some {\it  a priori} estimates for
the $\db $-equation, we need a variant of the
Bochner-Kodaira-Morrey-Kohn-H\"ormander formula (cf. Proposition 1.4)
as follows:

\proclaim{Proposition  3.1} Let  $\Omega$ be a domain in $\Bbb CP^n$ with
$C^2$-smooth  boundary
$b\Omega$. Let $\delta(x) = d(x, b\Omega)$ be the distance  function and
$\lambda=-(\delta)^t$
for some $t > 0 $.
  For any  $f \in  C^1_{(p,q)}(\overline \Omega)\cap \text{Dom}(\db^*)$,
where $0\le p\le n$ and $1\le q\le n$, we have
  $$\aligned   \int_{\Omega}  (-\lambda)|\bar\partial f|^2{} &  +
\int_{\Omega} (-\lambda)|\vartheta_{} f|^2{}+2\Re(\vartheta_{} f,
\overline{f\llcorner_{ (\db \lambda)_\#}}
 )
\\&= \int_{\Omega}  (-\lambda)\(| \overline\nabla f|^2+ \langle\Theta f ,
\bar f\rangle + \langle
(i\partial\db \lambda) f, \bar f\rangle\)
\endaligned
\tag 3.1.1 $$
  where
$\Theta$ and   $\overline\nabla $ are the same as in
  Proposition 1.4 and $\Re{( h)} $ is the real part of  $h$.
\endproclaim
\demo{Proof} The proof follows from the same calculation as in H\"ormander
\cite{H\"o2}.  There is no
boundary term since $t>0$.   Recall that $\vartheta = \db^*$.
A simple calculation shows that
$$
\vartheta_\phi f = \db^*_\phi f  = \db^* f  - f \llcorner_{( \db \phi
)_\# } \tag3.1.2
$$
holds for any $(p, q)$-form $u$ with $q \ge 1$. On the other hand,
by expanding the term $\delta^t  i\partial\db(-\log\delta^t )$  we
have
$$
\delta^t  i\partial\db(-\log\delta^t) = t  i\partial\db (- \delta^t) + t^2
\frac
{i\partial\delta\wedge\db\delta}{\delta^2}.
\tag3.1.3
$$

Substituting (3.1.2)-(3.1.3) into (1.33) of
Proposition 1.4 with $\phi = - t \log \delta$, we obtain (3.1.1). \qed
\enddemo

Using lemma 2.2 and Proposition 3.1, we have the following estimates:

\proclaim{Proposition  3.2} Let  $\Omega$ be a domain in $\Bbb CP^n$ with
$C^2$ boundary
$b\Omega$. Let $\delta$ be the distance function from $b\Omega$. Let $0<t
<1$ such that    $ \lambda=
-(\delta)^{t } $ satisfies condition (2.5).
  For any  $f \in L^2_{(p,q)}( \Omega)\cap \text{Dom}(\db )\cap
\text{Dom}(\db^*)$, where $0\le p\le n$ and $1\le
q\le n$, we have
  $$\aligned   &\int_{\Omega}(-\lambda)|\bar\partial f|^2      +
\int_{\Omega}(-\lambda)|\vartheta_{} f|^2{}  \\ &  \ge \hat C_t \(
\int_{\Omega}(-\lambda)|\overline\nabla f |^2{}
+\int_{\Omega}(-\lambda)|  f|^2+\int_{\Omega}(-\lambda)\frac{| 
 f\llcorner_{(\db\rho)_\#}|^2}{|\rho|^2} \),
\endaligned\tag  3.2
 $$
where $ \hat{C}_t$ is a constant number  independent of $f$.
 \endproclaim
\demo{Proof} Let $f\in C^1_{(p,q)}(\overline\Omega)\cap \text{Dom}(\db^*)$.
By (2.5), we have
$$  \langle
(i\partial\db \lambda) f, \bar f\rangle \ge C_1  \( \int_{\Omega}(-\lambda)|
f|^2+\int_{\Omega}(-\lambda)\frac{|  f\llcorner_{(\db\rho)_\#}|^2}{|\rho|^2}\).\tag 3.3$$
Also, we have for any $\epsilon>0$,
$$|2\Re (-\vartheta_{} f,    \overline{  f\llcorner_{(\db\lambda)_\#} }  )|\le
\frac t\epsilon
\int_{\Omega}(-\lambda)|\vartheta_{} f|^2 +t\epsilon
\int_{\Omega}(-\lambda)\frac{|  f\llcorner_{(\db\rho)_\#}|^2}{|\rho|^2} .\tag 3.4$$
Choosing $\epsilon$ sufficiently small,  (3.2) follows from
(3.1.1), (3.3) and Proposition 1.5.  This proves the
proposition for $C^1$-smooth forms.

 Since $C^1_{(p,q)}(\overline\Omega)\cap \text{Dom}(\db^*)$ is dense in
$\text{Dom}(\db )\cap
\text{Dom}(\db^*)$ in the graph norm (cf. Lemma 4.3.2 in \cite{CS}), we have
that $C^1_{(p,q)}(\overline\Omega)\cap
\text{Dom}(\db^*)$ is also dense in the graph norm with weights defined by
$(\int_{\Omega}(-\lambda)|\bar\partial f|^2     +
\int_{\Omega}(-\lambda)|\vartheta_{} f|^2+\int_{\Omega}(-\lambda)|
f|^2)^{\frac 12}$
from the Dominated Convergence Theorem. The proposition is proved by
  approximation of $C^1$-smooth forms. \qed
\enddemo

Recall that $u = \db^* N f $ is a solution of $\db u = f$ with the least
$L^2$-norm.
By the estimate above, we have the following:

\proclaim  {\bf Corollary 3.3} Let $\Omega\subset\subset\Bbb CP^n$
 and $t$ be the same as  in Proposition 3.2.
   Then
$$\p\bar\partial^* N_{(p,q)} f\p_{(t)} \le C \p f\p_{(t)}, \quad
f\in\text{Ker}(\db),
\ 2 \le q\le n.
\tag 3.5$$
$$\p\bar\partial  N_{(p,q)}  f\p_{(t)} \le C \p f\p_{(t)}, \quad
f\in\text{Ker}(\db^*),
\ 0\le q\le  n-1 .
 \tag 3.6$$
where $C$ depends only on $t$.
  \endproclaim
\demo{Proof}
 Notice that $N = \square^{-1}$ on the range of
$\square$. It follows that $N f \in \text{Dom}(\square)$. Hence,  by
(1.25.3),  we see that
 $\db^*Nf $ is in
$\text{Dom}(\db )\cap
\text{Dom}(\db^*)$.
If $q \ge 2$, for any $f\in L^2_{(p,q)}(\Omega)$,  $\db^*Nf$ is a
$(p,q-1)$-form. Thus for $f\in\text{Ker}(\db)$, one gets from (3.2)
that   $$\aligned
& \int_{\Omega}(-\lambda)| f|^2  =\int_{\Omega}(-\lambda)|\bar\partial
\db^*Nf|^2      +
\int_{\Omega}(-\lambda)|\db^*\db^*Nf|^2{} \\ &  \ge \hat C_t
\int_{\Omega}(-\lambda)|  \db^*Nf|^2,
\endaligned
 $$
where we used the following facts:  $f = (\db \db^* + \db^* \db) N f$,
$\db^* \db^* f = 0$ and
$\db N f = N \db f = 0$ with $f \in \text{ker}(\db)$.
This proves inequality (3.5) by choosing $ C = \frac{1}{\hat C_t}$.

Similarly, for any $0\le q\le n-1$, $\db Nf$
is in $\text{Dom}(\db )\cap
\text{Dom}(\db^*)$. Substituting $\bar\partial  N_{(p,q)}  f$ into (3.2),
for
$f\in\text{Ker}(\db^*)$, we have
  $$\aligned
& \int_{\Omega}(-\lambda)| f|^2 = \int_{\Omega}(-\lambda)|\bar\partial
\db Nf|^2      +
\int_{\Omega}(-\lambda)|\db^*\db Nf|^2     \\ &  \ge \hat{ C}_t
\int_{\Omega}(-\lambda)|  \db Nf|^2.
\endaligned
 $$
This proves the corollary  by choosing $ C = \frac{1}{\hat C_t}$.     \qed
\enddemo

Let $W^{s}_{(p, q)}(\Omega)$ be the Sobolev space  with $-\frac 12<s<\frac
12$ and let $\|\
\|_{s(\Omega)}$ denote its norm.
For any  $u$   in $\text{Dom}(\db )\cap
\text{Dom}(\db^*)$, we have $u\in W^1(\Omega,\text{loc}) $. This implies
that $u$ satisfies an elliptic system and
$u\in W^{s}(\Omega)$ for $-\frac 12 <s<\frac 12$    if and only if
$$\|u\|_{(-2s)}^2= \int_{\Omega}  \delta^{-2s} |u|^2 <\infty.\tag 3.7$$
For a proof of this, see Theorem C.4 in the Appendix in \cite{CS}.

  Using (3.7) or (2.9) and the estimates above, we are
ready to prove the Sobolev boundary regularity result for the (generalized)
Bergman  operator $  P: L^2_{(p,q)}(\Omega)
\to L^2_{(p,q)}(\Omega)\cap \text{ker}( \db) $. From Theorem 2.6,  we have
$P=I-\db^*N_{(p,q+1)}\db$.

 \proclaim  {\bf Theorem 3.4} Let $\Omega\subset\subset\Bbb CP^n$
be a pseudoconvex domain with  $C^2$-smooth boundary.  Let $0<t<t_0$ and
$t_0 = t_0(\Omega)$ be given by (0.7). Then  the Bergman projection  operator
 $P$ is bounded from  $W^{\frac t 2}_{(p,q)}(\Omega)$ to  $W^{\frac t
2}_{(p,q)}(\Omega)$, where $0 \le
p\le n$ and $0\le q\le n-1$.

\endproclaim
\demo{Proof}
By  Proposition 2.5,  we have that
    $\db^*N$ is bounded on Ker($\db$)   with
estimates
$$
\p\db^*  N_{(p,q)}  f\p_{(-t)} \le C ||| f|||_{(-t)'}  , \quad
f\in\text{Ker}(\db ),
\ 1\le q\le  n-1 .
 \tag 3.8
$$
 Let $P_t$ denote the Bergman projection with respect to the weighted space
$L^2( \delta^t)$.
For any $f, g\in L^2_{(p,0)}(\Omega)$ with $\db g=0$,
we have
$$
(Pf,g)=(f,g)=(\delta^{-t}  f,g)_{(t)}=(P_t\delta^{-t}
f,g)_{(t)}=(\delta^tP_t \delta^{-t}f ,g).
$$
This implies that
$$
\aligned P&=P^2=P\delta^tP_t \delta^{-t}
\\&=(I-\db^*N\db)\delta^tP_t \delta^{-t}\\&=\delta^t
P_t\delta^{-t}-\db^*N(\db(\delta^t)\wedge P_t\delta^{-t} ) \endaligned
\tag3.9
$$
since $\db P_t=0$.
For any $f\in L^2_{(p,q)}(\Omega)$,
$$\|\delta^tP_t \delta^{-t}f\|^2_{(-t)}\le \| P_t \delta^{-t}f\|^2
_{(t)}\le \|\delta^{-t}f\|^2 _{(t)}
=\|f\|_{(-t)}^2.\tag 3.10
$$
By (3.8), we have
$$
\aligned \p\db^*N(\db(\delta^t)\wedge P_t\delta^{-t}f )\p_{(-t)}^2&\le
C|||\db(\delta^t)\wedge P_t\delta^{-t}
f|||_{(-t)'}^2\\&\le C\p \delta^{\frac t2}   P_t\delta^{-t}
f\p^2=  C\p    P_t\delta^{-t}
f\p_{( t)}^2\\&\le C  \p  \delta^{-t}
f\p_{(t)}^2= C\p f\p_{(-t)}^2.   \endaligned
\tag 3.11
$$
From (3.9)-(3.11), we get
$$
\p Pf\p_{(-t)}\le C\p f\p_{(-t)}.\tag 3.12
$$
Note that $W^{\frac t2}\subset
L^2(  \delta^{-t})$ by (3.7) or (2.9). From (3.10), we get
$$
\p P   f\p_{(-t)}  \le C\p f\p_{(-t)}\le C_1\p f\p_{\frac t 2 }.\tag 3.13.1
$$
Using     (3.13.1) and (3.7) or (2.9), one obtains  that the Bergman projection
satisfies
$$\p Pf\p_{\frac t2}\le C_2 \p f\p_{\frac t2}. \tag3.13.2 $$
  Theorem 3.4 is proved. \qed
\enddemo

   Let us now study the Sobolev regularities for other operators: $N$, $\db
N$ and
$\db^* N$.

\proclaim  {\bf Theorem 3.5} Let $\Omega\subset\subset\Bbb CP^n$
be a pseudoconvex domain with  $C^2$-smooth boundary
$b\Omega$.  Let  $t_0 = t_0(\Omega)$ be the
order of plurisubharmonicity
for
the distance function  $\delta$ defined by (0.7). For any
$0<t<t_0$,   the $\db$-Neumann  operator
 $N$ is bounded from  $W^{\frac t 2}_{(p,q)}(\Omega)$ to  $W^{\frac t
2}_{(p,q)}(\Omega)$, where $0 \le
p\le n$ and
$0\le q \le n-1$. We also have the following estimates: for any $f\in
W^{s}_{(p,q)}(\Omega)$,
$$\p N f\p_{\frac t 2(\Omega)} \le 2 C^2 \p f\p_{\frac t 2(\Omega)},
 \tag 3.14$$
$$\p\bar\partial^* N f\p_{\frac t 2(\Omega)} \le C \p f\p_{\frac t 2(\Omega)},
 \tag 3.15$$
$$\p\bar\partial  N f\p_{\frac t 2(\Omega)} \le C \p f\p_{\frac t 2(\Omega)},
 \tag 3.16$$
where $C$ depends only on $t$.
  \endproclaim
\demo{Proof}
We already proved  that   the Bergman projection $P=I-\db^*N\db$ are
bounded on $W^{\frac
t2}.$
It is easy to verify  that, if $\Cal F^*$ is the adjoint map of $\Cal F$
with respect to the $L^2$-norm, then
$$
\split
& \|\Cal F u \p_{W^{\frac t2}(\Omega)}  = \sup_{v \in L^2} \{
\frac{( \Cal F u,  \bar v)_{L^2} }
{
\|  v \p_{W^{-\frac t2} }
}
\} \\
= & \sup_{v \in L^2} \{
\frac{(  u, \Cal F^* \bar  v)_{L^2} }
{
\|  v \p_{W^{-\frac t2}}
} \} \le \p \Cal F^* \p_{W^{-\frac t2}(\Omega) } \| u \|_{W^{\frac
t2}(\Omega)}, \\
\endsplit  \tag3.17
 $$
(compare with Lemma 1.3).

From the proof of  Theorem 3.4, we    have that the canonical solution
satisfies
 (3.8).
Since the  Bergman projection $  P$ is self-adjoint on $L^2(\Omega)$, by
Theorem 3.4 and
(3.17), we see that the
 Bergman projection $  P$ is bounded  on $W_{(p, q)}^{\pm \frac
t2}(\Omega).$

 Inequality (3.15) follows easily from Theorem 3.4 and the boundedness of the
Bergman projection   on $W^{\frac t2}$   since $\db^*Nf=\db^*NPf$.

 Let $\tilde P=\db^*N\db$ be the projection operator
into $\text{Ker}\db^*$. Then  $P =I-\db^*N\db=I-\tilde P$. It follows that
 $\db Nf=\db N\tilde Pf$.
 Since the
Bergman projection is self-adjoint, both  $P$ and $\tilde P$
satisfy
$$
\p Pf\p_{( t)}+\p \tilde Pf\p_{( t)}\le C_3\p f\p_{( t)}.\tag 3.18
$$
Thus from Corollary 3.3, we have that both $\db N$ and $\db^*N$ are bounded
with estimates
$$
\p \db N f\p_{( t)}=\p \db N \tilde Pf\p_{( t)}\le C_4\p \tilde Pf\p_{(
t)}\le C_4 C_3\p  f\p_{( t)},\ q\ge 0\tag
3.19
$$
$$
\p \db^* N f\p_{( t)}=\p \db^* N   Pf\p_{( t)}\le C_4 \p  Pf\p_{( t)}\le
C_4 C_3 \p  f\p_{( t)},\  q\ge 2.\tag 3.20
$$
In fact (3.20) also holds for $q=1$. This follows from the formula
$\db^* Nf=\db^*_tN_tf-P_t\db^*_tN_tf$ for any $f\in \text{Ker} (\db)$.
Thus for any $f\in L^2_{(p,1)}(\Omega)$, by (2.14) and (2.10.4) we have 
$$
\aligned
\p\db^* Nf\p_{( t)}&=\p\db^* N   Pf\p_{(
t)}\\&=\p\db^*_tN_t Pf-P_t\db^*_tN_t Pf\p_{(t)}\le C_5  \p P f\p_{(t)}\le C_5 C_3 \p
f\p_{(t)}.\endaligned\tag 3.21
$$
Notice that $(\db N)^* = N \db^* = \db^* N $ and
$(\db^* N)^* = N \db = \db N $. By (3.17)-(3.21), we obtain that
both $\db N$ and $\db^* N$ are bounded operators on $W_{(p, q)}^{\pm \frac
t2} (\Omega)    $. Thus, (3.15)-(3.16) are true, by choosing $C = \max\{C_3C_4, C_3C_5 \}$.

    It remains to verify (3.14). By an observation of Range,   we
   have  $$N=\db\db^*N^2+\db^*\db N^2=\db N \db^*N+\db^*N\db N. \tag3.22$$
It follows from  (3.15)-(3.16) and (3.22)
that
$$\p   N f\p_{( t)} \le 2 C^2 \p  f\p_{( t)},\  q\ge 0. \tag3.23$$
By (3.23) and (3.17), we conclude that $$\p   N f\p_{( -t)} \le 2 C^2 \p
f\p_{( -t)},\  q\ge 1, \tag3.24$$  since
$N$ is
self-adjoint. Thus (3.14) now follows from (3.24) and (3.7) or (2.9). \qed

\enddemo
Theorem 2 is a direct consequence of  Theorems 3.4-3.5.

\heading {\bf 4.  $\db $-closed extension from pseudoconvex boundaries in
$\Bbb CP^n$}\endheading

In this section we study the extension of $\db_b$-closed forms from the
boundary of a pseudoconvex domain
in $\Bbb CP^n$. This is equivalent to solving the $\db$-Cauchy problem on
pseudoconvex domains, which is the dual
of the $\db$-Neumann problem.

     Recall that if $\db_b v^{(0, 1)} = 0$ on $b\Omega$ and if $ w^{(0, 1)}$
is an arbitrary extension of $v^{(0, 1)}$ on $\Omega$,   then  we can
correct $w^{(0, 1)}$ to be
a $\db$-closed extension $\tilde v^{(0, 1)}$  on $\Omega$   by setting
$$
 \tilde v^{(0, 1)} = w^{(0, 1)} + {\star\bar\partial N_{(n,n-2)} {\star
(\db w^{(0, 1)})}} \text{       on  } \Omega, \tag4.0
$$
where $\star = \bar *:
L^2_{(p,q)}(\overline \Omega)$   is the Hodge star operator
defined by
$$\langle \phi,\bar\psi\rangle dV=\phi\wedge \star \bar\psi. \tag4.1
$$
Notice that $\star = \bar*$ satisfies $\star(\lambda u) = \bar{\lambda}
(\star u)$
for any complex number $\lambda  \in \Bbb C$.

    There are two issues in the application of the formula (4.0) above.
The first one is that we require
$\db w^{(0, 1)} \in L^2$   when we apply Theorem 3.5 to
the $(n, n-2)$-form $[\star{(\db w  ) }]$. Therefore, in
some results stated below,  we require that $ w \in W^1(\Omega)$ or
equivalently
$ v = w|_{b\Omega} \in W^{\frac 12  }(b\Omega)$, (cf. Proposition 4.3).

The second difficulty occurs when  $[\star{(\db w  ) }]$ is an $(n,
0)$-form, i.e., $n=2$.
It is easy to see that,  for a $(p, q)$-form $v$ on $M$  with $q < n-1$,
the $\db $-closed extension $ \tilde v$ on $\Bbb CP^n$ exists {\it only if
} $\db_b v = 0$. However, for all $(p,
n-1)$-forms
$v$ on $M$, the equation $\db_b v = 0$ always holds. There is another
necessary condition on a $(p, n-1)$-form
$v $ so that   its $\db$-closed extension on $\Omega$  exists, see (4.9)
(or equivalently  (4.3)) below. In order to  show
that (4.9) holds for any $(p, n-1)$-form
$v$ on the Levi-flat hypersurface $M$, we derive a new Liouville type
theorem for a pseudoconcave boundary
 $M = b\Omega$, see
Proposition 4.5 below.

To find a $\db$-closed extension $\tilde v$ on $\Bbb CP^n$
for $v$ with $\db_b v=0$ on $b\Omega$,  it suffices to solve
the inhomogeneous equation $\db u = f$ with compact support $\text{Supp}(u)
\subset\overline\Omega$.

\proclaim{\bf Proposition 4.1} Let $\Omega$ be a
pseudoconvex domain with $C^2$ boundary in $\Bbb CP^n$,
$n\ge  2$.
   For every $f\in L^2_{(p,q)}(\Bbb
CP^n)$, where $0\le p\le n$ and
$1\le q\le n-1$, with  $\bar\partial f=0$ in the
distribution sense in
$\Bbb CP^n$ and $f$ supported in $\overline\Omega$, one can
find $u\in L^2_{(p,q-1)}(\Bbb CP^n)$ such that
$\bar\partial u=f$ in the distribution sense in $\Bbb CP^n$
with
$u$ supported in $\overline\Omega$ and $$  \int\limits
_\Omega \mid u \mid^2 dV \leq C  \int\limits _\Omega \mid f \mid^2 dV  $$
for some $C>0$.
\endproclaim
\demo{Proof} From Theorem 2.6,  the
$\bar\partial$-Neumann operator of degree  $(n-p,n-q)$ in
$\Omega$, denoted by $N_{(n-p,n-q)}$, exists. Let $\star:
L^2_{(p,q)}( \Omega)$   be the Hodge star operator
defined as in (4.1) above.
   We
define
$$u=- {\star\bar\partial N_{(n-p,n-q)} {\star f}} =
- {\bar{*}\bar\partial N_{(n-p,n-q)} {\overline{* f}}},  \tag
 4.2$$
then $u\in L^2_{(p,q-1)} (\Omega) $.
By the proof of Corollary 3.3, we can verify that
$\star   u = \db N (\star f) \in \text{Dom} (\bar\partial^*) $,
because $N = \square^{-1}$ on $\text{Range}(\square)$ and  $\text{Range}(N)
\subset \text{Dom}(\square)$.
Since
$(\star u) \in \text{Dom}(\db^*)$, it follows from (1.25.2)  that
$u|_{T^{(0, 1)   }
(b \Omega   ) }= 0$.
Extending $u$ to $\Bbb CP^n$ by defining $u=0$ in $\Bbb CP^n
\setminus \Omega$, we obtain that $u$ has support in $\overline{ \Omega}$.
 It follows from a theorem of Kohn-Rossi \cite{KoR} (see also Theorem
9.1.2 in  Chen-Shaw \cite{CS})
that
$\bar\partial u=f$ in the distribution sense in $\Bbb CP^n$. \qed
\enddemo

In order to solve the $\db$-equation with
compact support
when $q= n$,   there is
another compatibility condition (see (4.3) below) and we
have the following result:

\proclaim {Proposition 4.2} Let $\Omega$ be a
pseudoconvex domain with $C^2$-smooth boundary  in $\Bbb CP^n$, $n\ge  2$.
For any
$f\in L^2_{(p,n)}(\Bbb CP^n)$, $0\le p\le n$,
such that
$f$ is supported in $\overline \Omega$ and
$$ \int_\Omega  f\wedge g=0\qquad \text{ for every}\quad
g\in L^2_{(n-p,0)} (\Omega) \cap\text{Ker}(\db),\tag
4.3$$ one can find $u\in L^2_{(p,n-1)}(\Bbb CP^n)$ such
that  $\bar\partial u=f$ in the distribution sense in $\Bbb
CP^n$ with $u$ supported in $\overline \Omega$ and
$$ \int\limits _\Omega \mid u \mid^2 dV \leq C\int\limits _\Omega \mid f
\mid^2 dV, $$
for some $C>0$.
\endproclaim
\demo{Proof} Using Theorem 2.6,    the
$\db$-Neumann operator
$N_{(p,0)}$ exists for any $0\le p\le n$ and
  we have
$$N_{(p,0)}=\db^*N^2_{(p,1)}\bar\partial. \tag 4.4$$ The Bergman projection
operator
$P_{(p,0)}$ is given by  $$\db^*\bar\partial
N_{(p,0)}=I-P_{(p,0)}.\tag 4.5$$ We define $u$ by
$$u=-\star {\bar\partial N_{(n-p,0)}{\star f}}.\tag
4.6$$
$$\aligned \db u&=(-1)^{p+n}
{\star\bar\partial^*\bar\partial N_{(n-p,0)}\star  f}\\&=
f-(-1)^{p+n} {\star P_{(n-p,0)}\star f}.\endaligned \tag
4.7$$

From (4.3), we get for any $g\in L^2_{(n-p,0)} (\Omega)
\cap\text{Ker}(\db)$,
$$ (\star f,  \bar{g})=(-1)^{p+n}\overline{ \int_\Omega  g\wedge
f}=0.  $$ Thus $P_{(n-p,0)}( {\star f)}=0$ and
$\bar\partial u=f$ in $\Omega$.
Using $\star u\in \text{Dom}(\db^*)$ and extending $u$ to
be zero outside
$\Omega$, then  $\bar\partial u=f$ in $\Bbb CP^n$ in the distribution
sense. \qed
\enddemo

   Let us now summarize the necessary and sufficient condition
on $f \in  W^{\frac 12}_{(p,q)}(b\Omega)$ to have a $\db $-closed
extension $F$ on $\Omega$.

\proclaim{Proposition 4.3} Let $\Omega\subset\subset\Bbb CP^n$ be a
pseudoconvex domain with  $C^2$ boundary
$M = b\Omega$.   Let
$f\in W^{\frac 12}_{(p,q)}(b\Omega)$, where $0\le p\le n$,  $ 1\le q\le n-1$.
We assume that
$$   \db_b f =0, \qquad  \text{if }1\le q<n-1, \tag 4.8 $$
and $$ \int_M f \wedge\psi=0 \quad \text{for every    }\psi\in L^2_{(n-p,
0)}(\overline \Omega)\cap\text{Ker}(\db),   \quad  \text{if }q=n-1. \tag 4.9 $$
  Then   there exists     $F\in L^2_{(p,q-1)}(\Omega)$ such
that
 $F=f$ on $ b\Omega $  and $\db F=0$ in $\Omega$.
\endproclaim
We remark that any $\psi\in L^2_{(n-p,
0)}(  \Omega)\cap\text{Ker}(\db)$ has   boundary values on
$b\Omega$ in the $W^{-\frac 12}(b\Omega)$ space
(see e.g. Lemma 2.1 in Michel-Shaw \cite{MS2}). Thus  the pairing (4.9) is
well-defined in the sense of currents.
\demo{Proof}
  Since $f\in W^{\frac 12}_{(p, q)}(b\Omega)$ is a form, one can extend $f
= \Sigma_{|I| = p, |J| = q} f_{I, \bar J}$
componentwise to $\Omega$ such that each component
$ f_{I, \bar J}$ is in $W^1(\Omega)$.  For
detailed construction of such an extension, see e.g.  Lemma 9.3.3 in
\cite{CS}. Let
$\tilde f$ be an arbitrary
extension
 of $f$ with $ \tilde f \in W^1_{(p, q)}(\Omega)$

We first assume that  $q+1<n $.      From (4.8), we can require  that
$f_1=\db_b
\tilde f=0$ in $M$. If we extend $f_1$ to be zero outside $\Omega$, we get
 $\db f_1=0$ in $\Bbb CP^n$ in the distribution sense.
  we set $v_0=- {\star\bar\partial N_{(n-p,n-q-1)} {\star f_1}}$. From
Proposition 4.1 and its proof, we have   $v_0\in
L^2_{(p,q)}(\Omega)$. We set $v_0=0$ outside $\Omega$. Then  $\db v_0=f_1$
in the sense of distributions on  $\Bbb CP^n$   with
 $ \text{supp}( v_0)
\subset\overline \Omega$.

Setting $F=\tilde f-v_0$ in $\Omega$, we have $F=f$ on $b\Omega$ and $\db
F=0$ in $\Omega$. This proves the
proposition for $q<n-1$.

When $q=n-1$ and $f$ satisfies (4.9), we let $f_1 = \db \tilde f$.  Then for any $\psi\in
C^1_{(n,0)}(\overline\Omega)$
with
$\db \psi=0$, we have
$$\int_\Omega f_1\wedge \psi=\int_\Omega \db \tilde f \wedge
\psi=\int_{b\Omega}f\wedge \psi=0.\tag 4.10
$$
 When $ \psi \in L^2_{(n,0)}(\Omega) \cap \text{ker}( \db) $, (4.10)
follows from approximating
$\psi$ by smooth forms. Let us now apply Proposition 4.2 to $f_1$,  since
$f_1$ satisfies (4.3) by (4.10).
Setting $v_0 =- \star\bar\partial N_{(n-p,n-q-1)} {\star f_1}$ and $F
=\tilde f-v_0 $ as above, by
the proof of Proposition 4.2, we can conclude the proof for the case of $q=
n-1$ as well.
\qed
\enddemo

Recall that the domain $\Omega_+$ is pseudoconcave if and only if
its complement $\Omega_-=  \Bbb CP^n \setminus \overline\Omega_+ $ is
pseudoconvex.
 In order to show that the condition (4.9) automatically holds
on the Levi-flat boundary $M = b\Omega_{\pm}$, we observe that
$M$ is the boundary of a {\it pseudoconcave} domain as well. Indeed,
if $\Omega_+\cup \Omega_- =  \Bbb CP^n \setminus M $ and if M is Levi-flat,
then $\Omega_{\pm}$ is both pseudoconvex and  pseudoconcave.

 Notice that, by a theorem of Chern-Lashof,
any compact domain $\bar\Omega $ with $C^2$-smooth boundary $b\Omega$ in $
\Bbb C^n$
must have a point $Q_0$ in $b\Omega$ with strictly positive principle
curvature.
To see this fact, one considers the Gauss map from $b\Omega $ to the unit
sphere $S^{2n-1} \subset \Bbb C^n$.
Therefore, this theorem of Chern-Lashof implies that there is {\it no}
compact pseudoconcave
domain $\overline{\Omega}$ in $\Bbb C^n$. However, the complex projective
space $ \Bbb  CP^n $ contains
a lot of  pseudoconcave
domains.

The classical Liouville Theorem states that
if $f: \Bbb C^n \to \Bbb C$ is a bounded holomorphic function on
$ \Bbb C^n$ then $f$ must be constant. We will extend the classical
Liouville Theorem of functions  to the $(p, q)$-forms on pseudoconcave
domains $\Omega^+$ as
follows.

First  we study the one-sided $\db$-closed extension from a
pseudoconcave domain $\Omega_+$ to the whole space $\Bbb CP^n$
by using (4.8).

\proclaim{Proposition 4.4} Let $\Omega_+ \subset\Bbb CP^n$ be a
pseudoconcave domain with $C^2$-smooth
boundary.
Suppose that
$f\in W^1_{(p,q)}(\overline{\Omega}_+)$  with $\db f = 0$, where
$0\le p\le n$ and  $ 1\le q< n-1$.  Then there exists a $\db$-closed form
$F\in L^2_{(p,q)}(\Bbb CP^n)$ such that
$F$ is a $\db$-closed extension of $f$.
\endproclaim
\demo{Proof} Let $M = b\Omega_+$ be the $C^2$-smooth pseudoconcave boundary
in $\Bbb CP^n$. We consider the so-called Fermi map (cf. [Cha])
along $M$ as follows:
$$
\split
      h: \quad  & M \times \Bbb R \to \Bbb CP^n \\
          & (Q_0, t) \to Exp_{Q_0}(t \nabla \rho),
\endsplit \tag4.11
$$
where $Exp$ is the exponential map of $\Bbb CP^n$ and $\rho(x) =
- d(x, b\Omega_-)$ for $x\in \Omega_-$. It is well-known
that, for sufficiently small $\epsilon >0$, the Fermi-map
$h|_{M \times [-\epsilon, \epsilon]}$ is an embedding map.

Let us now choose an extension $w \in L^2_{(p, q)}(\Bbb CP^n)$
of $f$ as follows: for $0 \ge  t \ge -\epsilon$, let
$$
w(h(Q_0, t)) = f(h(Q_0, -t)). \tag4.12
$$
Thus, $w$ is well-defined on $U_\epsilon(M)$. We further extend
$w $ to $ [\Omega_- \setminus U_{\frac{\epsilon}{2}}(M)] $ so that
$w|_{\overline{\Omega}_+  } = f$ and $w \in W^1_{(p,q)}(\Bbb CP^n)$.

Recall that $\Omega_-$ is pseudoconvex by the assumption.
Applying Proposition 4.1 to the equation
$$
\db u = \db w \quad \text{  on } \Omega_-,
$$
we get
a solution $u$ with  support in $\overline{\Omega}_-$ by choosing
$$
  u = - {\star\bar\partial N_{(n-p,n-q-1)} {(\star  \db w)}} \quad \quad
\text{    on } \Omega_-. \tag4.13
$$
Inspired by
(4.0), we let
$$
   F_- = w - u = w +  {\star\bar\partial N_{(n-p,n-q-1)} {(\star  \db w)}}
\quad\quad     \text{    on } \Omega_- , \tag4.14
$$
where
$u $ is given by (4.13). A direct computation shows that
$$
\db u = \db w - (-1)^{p + q} \star P_{(n-p, n-q)} \star (\db w), \tag4.15
$$
where $P_{(n-p, n-q)}: L^2_{(p, q)}(\Omega_-) \to \text{Ker}(\db )$ is the
Bergman  projection.

   When $w$ is a $(p, q)$-form with $q < n-1$, $\db w \in\text{Ker}(\db)$
and hence $P(\star(\db w)) = 0$. By (4.15) and $P(\star(\db w)) = 0$, we
have  that $\db u = \db w$
and $\db F_- = 0$.

By the proof of Proposition 4.1, we know that $u$ has compact support in
$\overline{\Omega}_-$.
 Define
$$F=\cases &f,\quad x \in
\bar{\Omega}_+,\\&F_- ,\quad x\in  \Omega_-.\endcases \tag 4.16$$
Then $F\in L^2_{(p,q)}(\Bbb CP^n)$, $F=f$ on $\bar{\Omega}_+$ and $\db F=0$ in $\Bbb
CP^n$ in the distribution sense.
This proves the Proposition. \qed
\enddemo

In order to show that (4.9) holds, it is sufficient to derive a new
Liouville type Theorem, i.e., to show that
$ L^2_{(n-p,0)}(\Omega) \cap \text{Ker}( \db) = 0 $
for a pseudoconcave domain     and  $n-p >0$.

\proclaim{Proposition  4.5} Let $\Omega_+\subset\subset\Bbb CP^n$ be a
pseudoconcave domain with  $C^2$-smooth boundary
$b\Omega_+ $.   Then $L^2_{(p,0)}(
\Omega_+)\cap \text{Ker}(\db)=\{0\}$  for every
$0< p\le n$;  and  $L^2_{(0,0)}(
\Omega_+)\cap \text{Ker}(\db)=\Bbb C$.
\endproclaim
\demo{Proof}
 We will extend $f$ to be  a $\db$-closed form in $\Bbb
CP^n$ in a similar way as in the proof of Proposition 4.4. If $f\in
L^2_{(p,0)}( \Omega_+)\cap \text{Ker}(\db)$,  then
  there
exists a   $\db$-closed
extension $F \in L^2 $  on the whole $\Bbb CP^n$. To see this, we extend
$f$ to be
$\tilde f$ in
$\Bbb CP^n$ such that
$f_1\equiv\db \tilde f\in L^2_{(p,1)}(\Bbb CP^n)$ with support
$f_1\subset\overline{\Omega}_-$, where $\Omega_-=\Bbb
CP^n\setminus \overline\Omega_+$.

 Such an $L^2$-extension $\tilde f$ can be carried out,  because $b\Omega$
is a
$C^2$-smooth real hypersurface. Here is the detail for the construction of
$\tilde f$ and $f_1$.
Let $h: b\Omega \times [-\epsilon, \epsilon]
\to \Bbb CP^n $ be the Fermi-map given by (4.11). For sufficiently
small $\epsilon$, $h$ is an embedding map which gives rise to the so-called
Fermi coordinate along $b\Omega$. We now define $ \tilde f(h(Q, t)) = f
(h(Q, -t))$
for $h(Q, t) \in  \Omega_-$ and $|t| < \epsilon/2 $. Outside the annuli
neighborhood
$h( b\Omega \times [-\epsilon, \epsilon]     ) $, we require
$\tilde f$ to be $C^2$-smooth. Therefore, the $(p, 1)$-form $f_1\equiv\db
\tilde f \in L^2_{(p,1)}(\Bbb CP^n)$ with support
$f_1\subset\overline{\Omega}_-$.

  Using Propositions 4.1-4.4, we see that   there exists
a   solution $v_0\in L^2_{(
p,0)}(\Omega_-)$ such that
$\db v_0=f_1$ in $\Bbb CP^n$ if we extend $v_0$ to be zero on $\Omega_+$.
Thus,  the extension  form $F =\tilde f -v_0$ is   $\db$-closed in the
whole space
$\Bbb CP^n$. We observe that  \newline $L^2_{( p, 0)}(\Bbb CP^n)\cap
\text{Ker}(\db)
 =\{0\}$ for all $0< p
\le  n$;
and $L^2_{( 0, 0)}(\Bbb CP^n) \cap \text{Ker}(\db) = \Bbb C$. Hence, $f =
0$ when $p> 0$,  and $f$ must be   a
constant
when it is a $L^2$-holomorphic function on $\Omega_+$.
This proves Proposition  4.5. \qed
\enddemo

  Some related results similar to Propositions 4.4-4.5 were obtained by
Henkin-Iordan [HI] earlier via a different argument.
Let us now apply our results above to a $C^2$-smooth
Levi-flat hypersurface $M$ if it exists.

\proclaim{Corollary 4.6} Let $M\subset\Bbb CP^n$ be a   $C^2$-smooth
Levi-flat hypersurface.
For any
$f\in W^{\frac 12}_{(p,n-1)}(M)$, where $0\le  p< n$,
 there exists a $\db$-closed form $F\in L^2_{(p,n-1)}(\Bbb CP^n)$ such that
$F$ is a $\db$-closed extension of $f$. \endproclaim
\demo{Proof}
Let $M$ be a compact connected $C^2$-smooth real hypersurface in $\Bbb
CP^n$. It is well-known
that $f \in  W^{\frac 12}_{(p,q)}(M)$ if and only if there exists an
extension $\tilde f$ of $f$ with
$\tilde f \in  W^{1}_{(p,q)}(\Bbb CP^n)$.

   Since $\Bbb CP^n$ is simply-connected, $\Bbb CP^n\setminus M = \Omega_+
\cup \Omega_-$ has exactly two connected components: $\Omega_+$ and
$\Omega_-$. Since $\Omega_{\pm}$ is pseudoconvex, we let
$$
 F_{\pm} = \tilde f +  \star[\bar\partial N_{(n-p,0)}|_{\Omega_\pm}]
{(\star  \db \tilde f)} \quad  \quad \text{        on } \Omega_{\pm}
\tag4.17
$$
and $ v_{\pm} = - \star\bar\partial N_{(n-p,0)} {(\star  \db \tilde f)}$.
By Theorem 2.6, we have
$F_\pm \in L^2_{(p,n-1)}(\Bbb CP^n)$.

   Since $M$ is Levi-flat, we see that both $\Omega_+$ and $\Omega_-$ are
pseudoconcave as well.
By Proposition 4.5, the form $f_1 = \db \tilde f $ satisfies (4.3) and
(4.9) on each of
$\Omega_\pm$. By the proofs of Propositions 4.1-4.5, we conclude that $F$
is the $\db $-closed extension of
$f$ on the whole $\Bbb CP^n$. \qed
\enddemo

 We are now ready to study the inhomogeneous equation (0.4).

\proclaim{Corollary  4.7} Let $M\subset\Bbb CP^n$ be a $C^2$-smooth
Levi-flat hypersurface. Let
$f\in W^{\frac 12}_{(p,q)}(M)$,  where $0\le p\le  n $, $1\le q\le n-1$ and
$p\neq q$. We assume that $f$ is
$\db_b$-closed in
$M$ if $q<n-1$ and $p<n-1$  if $q=n-1$.  Then there exists
$u\in W^{\frac 12 }_{(p,q-1)}(M)$ such that
 $\db_b u=f$ in $ M$.
\endproclaim

\demo{Proof}
From Proposition 4.4 (for $q<n-1$) and Corollary 4.6 (for $q=n-1$), there
exists a $\db$-closed  $F\in L^2_{(p,q)}(\Bbb
CP^n)$ such that $F=f$ on $M$.

 From our assumption $p\neq q$, we have that  $H^{p,q}(\Bbb CP^n)=\{0\}$.
This implies that
there exists an inverse
$\Bbb G$,  the Green operator, for
$\square$ on $\Bbb CP^n$ such that
$\square \Bbb G =I$ on $L^2_{(p,q)}(\Bbb CP^n)$
  and the
   the Hodge decomposition
theorem holds for $(p,q)$-forms on $\Bbb CP^n$. For any
$\db$-closed $(p,q)$-form $F$,  we have
$$
F=\db\db^*\Bbb GF+\db^*\db \Bbb GF=\db\db^*\Bbb GF+\db^*  \Bbb G\db
F=\db\db^*\Bbb GF . \tag4.18
 $$
 Thus $F=\db \tilde u$ where
$ \tilde u =\db^*\Bbb GF$ is the canonical solution   on $\Bbb CP^n$.
Using the interior regularity for $\db$, we see that $
\tilde u \in W^1(\Bbb CP^n)$.
Restricting $\tilde u$ to $M$ and denoting the
restriction by $u$, we have $u\in W^{\frac 12}(M)$ and $\db_bu=f$.
\enddemo

\heading {\bf 5. Proof of  Theorem 1}\endheading

We will use results of Section 4 to prove Theorem 1.
Let $M$ be a $C^{2,\alpha}$ Levi-flat hypersurface in $\Bbb CP^n$, $n\ge 2$
and let   $\rho=-\delta$ be the
signed distance function for $M$.
We first discuss   solvability with regularity for the   $\db_b$ equation
in the Sobolev
spaces on $M$.

As we pointed out  in Section 0, the proof of Theorem 1 is reduced to
a problem of finding a continuous solution $u $ of
$$
i\partial_b\db_b u = f_b = f|_{  [T(M)]_{\Bbb R}
\cap J [T(M)]_{\Bbb R}  } \quad \text{ on } M^{2n-1} \tag5.1$$
on $M$
under the assumption that
$
f = d v
$
is an exact form and $f_b$ is a real-valued H\"older continuous (1,
1)-form. When
$f_b$ is a  (1, 1)-form and
$v^{(0, 1)}$ is part of $v$, we have
$
   \db_b v^{(0, 1)} = 0 \text{ and } \partial_b v^{(1, 0)} = 0.
$

    There is the corresponding
classical Lelong equation
$$
  i \partial\db \tilde u = \tilde f \quad \text{ on } \Bbb CP^n \tag5.2
$$
where
$
   \tilde f = d \tilde v
$
is an exact real-valued $(1,1)$-form on $\Bbb CP^n$. Since $\tilde f$ is of
type $(1,1)$, we have
$
   \db \tilde v^{(0, 1)} = 0 \text{ and } \partial \tilde v^{(1, 0)} = 0
$
where $\tilde v = \tilde v^{(0, 1)} + \tilde v^{(1, 0)}$ on
$\Bbb CP^n$.
The $\partial\db$-exact
Poincar\'e Lemma states that if  $\Bbb G_{(0, 1)}$ is the Green's operator
of the Laplace operator
 $\square_{(0, 1)}$ acting on $(0, 1)$-forms of $\Bbb CP^n$, then
$
  \tilde u = 2 \Re \{ \db^* \Bbb G_{(0, 1)}  \tilde v^{(0, 1)} \}$ on $\Bbb
CP^n
$
is a solution of (5.2), where $\Re \{ \phi\}$ denotes the real
part of $\phi$ (cf. [Zh]).

Hence, to solve our equation (5.1) it is sufficient
to extend the
$\db_b$-closed 1-form $v^{(0, 1)}$ of $M$ to be a global $\db$-closed form $
 \tilde v^{(0, 1)}$ on the whole space $\Bbb CP^n$.
Recall that  if  $w^{(0, 1)}_{\pm}$ is an {\it arbitrary} extension of
$v^{(0, 1)}$ from $M = b\Omega_{\pm}$ to $\Omega_\pm$,
the $\db$-closed extension $\tilde v^{(0, 1)}_{\pm}$ of $v^{(0, 1)}$ is
given by
$$
  \tilde v^{(0, 1)}_{\pm} = w^{(0, 1)}_\pm  +
 [ \bar{*}(\db   N_{(n, n-2)}|_{\Omega_\pm})]\bar * (\db w^{(0, 1)}_{\pm}).
$$

We would also like to  explain why the  $W^s(\Omega)$-regularity result of
Theorem 2 is good enough
for the proof of Theorem 1. When the Levi-flat hypersurface $M$ is $C^{2,
\alpha}$-smooth,
we can choose $w = \theta $ to  be the connection form for the complex line
bundle
of equidistant hypersurfaces in a neighborhood $U_\epsilon(M)$ of  $M$. It
was shown in Section 1 that both
$w$ and $dw$ (i.e., $\theta$ and the curvature $\Theta^N = d\theta$) are
$C^{0, \alpha}$-smooth
in  $U_\epsilon(M)$. If Theorem 2 holds,  the
$\db $-closed extension $\tilde v^{(0,1)}$ of $\theta^{(0,1)}$ is in
$W^s_{(0, 1)}(\Bbb CP^n)$.
Using  elliptic theory on $\Bbb CP^n$, we obtain that the solution
$\tilde h = \tilde u$ to (5.2) is in $W^{1+s}_{(0, 0)}(\Bbb CP^n)$, see
[GT]. Applying
the trace theorem in Sobolev spaces   to $\tilde h$, we conclude
that $h = \tilde h|_M$ is in $W^{\frac 12 +s}(M)$. By the classical
elliptic theory on
the holomorphic leaves of $M$, we  already know that
$h $ is $C^\infty$ in each holomorphic leaf of $M$. With some
extra efforts, we can show that $h $ is H\"older continuous in remaining
directions of
$[T(M)]_{\Bbb R}$, whenever $h \in W^{\frac 12 +s}(M)$ and
$ i \partial_b {\db}_b h = \Theta_b \in C^{0, \alpha}(M) $ for some $s,
\alpha> 0$,
see Lemmas 5.1-5.2  below.

\proclaim{Lemma 5.1}  Let $M\subset\Bbb CP^n$ be a $C^2$-smooth Levi-flat
hypersurface and let
$f\in W^{\frac 12}_{(p,q)}(M)$,  where $0\le p\le  n $, $1\le q\le n-1$ and
$p\neq q$. Suppose that  $\epsilon_0= \frac 12  \min\{t_0(\Omega_+),
t_0(\Omega_-)\}$
and $t_0(\Omega)$ is the order of  plurisubharmonicity of $\Omega$ given by
Definition 0.1.
 We further assume that   $f\in W^{\frac 12+\epsilon}_{(p,q)}(M)$ is
$\db_b$-closed in
$M$ if $q<n-1$ and $p<n-1$  if $q=n-1$.
   Then   for any  $0\le \epsilon<\epsilon_0$,
there exists
$u\in W^{\frac 12+\epsilon}_{(p,q-1)}(M)$ satisfying
 $\db_b u=f $ in $ M.$

\endproclaim
\demo{Proof} Let $\Bbb CP^n\setminus M=\Omega_+\cup \Omega_-$.
Let $t_0^+$ and $t_0^-$ be the orders of plurisubharmonicity  for the
distance functions associated  with $\Omega_+$ and
$\Omega_-$ respectively defined in Definition 0.1.

We use the same notion as in Section 4. We can construct
a $\db$-closed form $F\in
W^{\epsilon}_{(p,q)}(\Bbb CP^n)$ such that $F|_M = f$. This is proved by
Theorem  3.5 and (4.14)-(4.18) as follows.

Since $f\in W^{\frac 12+\epsilon}_{(p,q)}(M)$,
there exists an extension  $\tilde f \in W^{1+\epsilon}(\Bbb CP^n) $ of $f$
by the Trace Theorem (cf. Appendix of [CS]). As in the proof of Corollary
4.6,  we let
$$
 F_{\pm} = \tilde f  +  \star[\bar\partial N_{(n-p,n-q-1)}|_{\Omega_\pm}]
{(\star  \db \tilde f)}  \quad \quad  \text{    on } \Omega_{\pm}.
$$
This
implies that  $F\in W^{\epsilon}(\Bbb CP^n)$ by Theorem 3.5.

Next, we solve the equation $\db \tilde u = F$ on the whole manifold $\Bbb
CP^n$.
 For the same reason as  in the proof of Corollary  4.7, we conclude that
$ \tilde u \in
W^{1+\epsilon}$ with $\db  \tilde u =F$. This gives  that $u= \tilde u|_M$
is in
$W^{\frac 12+\epsilon}_{(p,q-1)}(M)$ by the Trace Theorem again. \qed
\enddemo

We need two more preliminary results for the proof of Theorem 1.
We denote by  $\Cal N_{{ \widetilde{\nabla\rho}  }}^{1,0}$  the complex
line bundle
and its associated curvature form by $\tilde\Theta^N$ as in Section 1.
From Propositions 1.1-1.2,  $\Theta_b$ is a
$(1,1)$-form on the Levi-flat hypersurface $M$. As in Section 1, we let
$\beta$ be the
real 1-form such that
$$ \tilde\Theta^N=-d\theta=\sqrt{-1}d\beta  \quad \text{ on } U_\epsilon(M), $$
where
$$\beta(\cdot)=\text{Hess}(\rho)(\cdot,J(\nabla\rho)).\tag 5.3$$  Write $$
\beta=\beta^{1,0}+\beta^{0,1}$$ where
$\beta^{1,0}$ and
$\beta^{0,1}$ are the (1,0) and (0,1) components of $\beta$.

\proclaim{Lemma 5.2} Let $M$ be a compact $C^{2,\alpha}$ Levi-flat
hypersurface in $\Bbb CP^n$ and
$\beta^{0,1}$ be the $(0,1)$-form defined in (5.3).  Then
there exists a function $u\in C^{ \epsilon}(M)$  such
that
 $$\db_b u=\beta^{0,1}
\quad
\text{in } M \tag 5.4$$
for sufficiently small $\epsilon < \frac 14 \min\{ \alpha,
t_0(\Omega_\pm)\}$, where
$t_0(\Omega_\pm) $ is the plurisubharmonicity of $\Omega_\pm$ and $\Omega_+
\cup \Omega_- = \Bbb CP^n$.
\endproclaim
\demo{Proof} Since $\rho$ is of class $C^{2,\alpha}$,
the 1-form $\beta$ is obtained from the Hessian of $\rho$, hence, is in
$C^\alpha$.
The Levi-flat hypersurface $M$ is locally  foliated by complex manifolds
$\{ \Sigma_t\}$, where $V = \cup_{|t| < \mu} \Sigma_t$ is an open subset of
$M^{2n-1}$. Since each leaf $\Sigma_t$ is a complex hypersurface, it follows
from the Chern formula that  the curvature tensor    $\Theta_b = \tilde
\Theta^N|_{\Sigma_t}$
 is a
$(1,1)$-form (see Proposition 1.2).
  Thus from (1.11) and type consideration, we get
for $n>2$,
$$\db_b\beta^{0,1}=0\quad \text{in } M. \tag 5.5$$
When $n =2$, $\beta^{0,1}$ satisfies the compatibility condition
(4.9) by Proposition 4.5.

We claim that there exists a   solution $u$ of (5.4) such that $u\in
W^{\frac 12+\epsilon}(M)$ for some
$\epsilon>0$.   If $\alpha\ge \frac 12$, we can use Lemma 5.1 directly. If
$\alpha<\frac 12$,
we note  that the Cartan-Chern-Gauss structure equation holds in a tubular
neighborhood $U = U_\epsilon(M)$ of $M$.
Thus,
$\beta^{0,1}\in C^\alpha(U)\subset W^\alpha(U)$.  Using the last assertion
of  Proposition
1.1,    $\db \beta^{0,1}$ is in $C^\alpha(U)$ since it is the
$(0,2)$-component of the
curvature form $\tilde\Theta^N$, which is in $C^\alpha$.  Thus   $
\beta^{0,1}$ has an extension to $\Bbb CP^n $
with
$\db\beta^{0,1}$  in $W^{\alpha}_{(0,2)}(\Omega_+)$ and
$W^{\alpha}_{(0,2)}(\Omega_-)$ respectively.

 Repeating the same arguments as in the proof of Proposition 4.4 for
$q<n-1$  and Corollary 4.6 for $q=n-1$,   we obtain a
$\db$-closed extension
$F \in L^2$ of
$\beta^{0,1}$ on the whole space $\Bbb CP^n$ with
  $F = \beta^{0,1}$  on $M$. Also from the boundary regularity of the
$\db$-Neumann operator proved in Theorem 3.5,  the
extension
$F$ is in $W^\epsilon(\Bbb CP^n)$ with
$\epsilon<\min\{\alpha,\epsilon_0\}$, where $\epsilon_0$ is  defined in
Lemma 5.1.

Using
the same arguments as in the proof of
Corollary 4.7, we can prove that  there exists  $\tilde u \in W^{1 +
\epsilon } (\Bbb CP^n)$ satisfying $\db \tilde u = F$. Therefore, if $u = \tilde u|_
M$, then
we have
$u\in
W^{\frac 12+\epsilon}(M)$   satisfying (5.4).

    It remains to  show that $u$ is H\"older continuous in $M$. This follows
from the following version of the Sobolev embedding
theorem. We note  that, on each leaf $\Sigma_t$,
$u$ satisfies  an elliptic equation. Thus we already have that $u$ is
smooth on each leaf $\Sigma_t$,  because $(\partial \db u)|_{\Sigma_t} =
\Theta_b$
and $\Theta_b$ is analytic on each holomorphic leaf $\Sigma_t$.
It remains only to show that
$u$ is H\"older continuous in the transversal
direction $\frac{\partial}{\partial t}$. To do this we need
to parametrize our hypersurface $M$ locally.

  Let $(z',g(z',t))$ denote the leaf
$L_t$ where $g(z',t)$ is holomorphic in $z'\in \Bbb B^{n-1} \subset \Bbb
C^{n-1}$ and
$C^2$-smooth  in $t$.
We can parametrize $M$ locally as a graph of a function $\eta + g$, by setting
$$\Psi(z',t)=(z',\eta(t)+g(z',t)),
$$ where
$z'\in \Bbb C^n$,
 $0\le |t|<\mu$ and $\eta(t)$ is a   $C^{1,\alpha}$ function in $t$ with
$\eta(0)=0$ and
$\eta'(0)=1$.   Clearly, $\Psi: \Bbb B^{n-1} \times (-\mu, \mu)
\to M$ is a local coordinate map of $M$, where $\Bbb B^{n-1} \subset
\Bbb C^{n-1}$ is an open set of $\Bbb C^{n-1}$.  Using a result of
Barrett-Fornaess (see \cite{BaF}), the foliation is actually $C^{2,\alpha}$.

It is easy to see that the push forwards
$\bar L_i=\Psi_*(\frac{\partial\  }{\partial \bar z_i})$ of $\frac{\partial\
}{\partial \bar
z_i}$,  $i=1,\cdots, n-1$,
  are the tangential Cauchy-Riemann equations for $M$. Let
$$
T=\Psi_*(\frac{\partial \ }{\partial t
}).\tag 5.6$$ Since
$$\[\frac
{\partial \ }{\partial t } ,\frac{\partial\  }{\partial \bar z_i} \]=0,\
\[\frac
{\partial \ }{\partial t } ,\frac{\partial\  }{\partial   z_i} \]=0,\tag 5.7$$
we have
$$\[
T ,\bar L_i \]=0,\ \[T ,L_i\]=0.$$
Thus the tangential Cauchy-Riemann equations $\db_b$ are just the
Cauchy-Riemann equations on $\Sigma_t$  and they commute
with
$T$.  When $u$ is restricted to each leaf $\Sigma_t$, $u$ satisfies an
elliptic system
in coordinates $z'$ and
$$\partial_{z'}\db_{z'} u= \Theta_b \quad \quad  \text{  on } \Sigma_t,\tag
5.8$$
where  $\Theta_b = \tilde{\Theta}^N|_{[T(M)]_{\Bbb R}
\cap J [T(M)]_{\Bbb R} }= \tilde{\Theta}^N|_{\Sigma_t}  $ is a  $(1,
1)$-form and  $J$ is the complex structure of
$\Bbb CP^n$, see Proposition 1.2.

From the classic Schauder theorem (cf. [GT]) for elliptic systems on
$\Bbb C^{n-1}$,  we get that $u\in
C^{2,\alpha}(\Sigma_t)$ for each $t$ since
$\tilde\Theta^N$ is in $C^\alpha$.  Furthermore, we have that there
exists a constant $C_1$   independent of
$t$    such that
$$|u|_{C^{2,\alpha}(\Sigma_t)}\le C_1(|\tilde\Theta^N|_{C^{
\alpha}(\Sigma_t)}+\|u\|_{L^2 (\Sigma_t)}),\tag 5.9$$
where $C_1$ depends on $V$ and $\Psi$, but $C_1$ is independent of $t$,
(because the local foliation is $C^2$-smooth
and the equation (5.8) is uniformly elliptic on $\Sigma_t \subset V$
independent of $t$).

Recall that $V = \cup_{|t| < \mu} \Sigma_t \subset M$. From the Sobolev
trace theorem,  a function $u\in W^{\frac 12+\epsilon}(M) $ has
$L^2$-trace on each leaf. Therefore,
 there exists $C_2>0$  independent of $t$ such that
$$\|u\|_{L^2 (\Sigma_t)}\le C_2 \|u\|_{W^{\frac 12+\epsilon} (M)}.
\tag 5.10$$
Combining (5.9) and (5.10), we get
$$   |u|_{L^\infty(V)}\le \underset{|t| < \mu}\to\sup
|u|_{C^{2,\alpha}(\Sigma_t)}\le C_3.\tag 5.11 $$
Thus we have  already proved that  $u$ is bounded.

It remains to  show that
$u$ is H\"older continuous in $t$
as well. To do this we differentiate the equation (5.8) in $t$ with order
$0 < \epsilon < 1$.

Let  $D_{t,h}^\epsilon $   denote the
finite difference
$$D_{t,h}^\epsilon u= \frac{  u(z',t+h)- u(z',t)  }{ |h|^\epsilon} $$
 and let
$|\ |_{  s}$ be the Besov norm  given by
$$|u|_s=\underset{0<|h|<\eta_0}\to\sup \frac{ \|u(x+h)-u(x)\|}{|h|^s}.$$
The Besov norm is weaker than  the Sobolev norm $\|u\|_s$.  It is easier to
use the Besov norm than the Sobolev
norm (see e.g. H\"ormander \cite{H\"o4})
 and we have   $$\|u\|_{s'}\le |u|_s\le \|u\|_s\tag 5.12$$
for any $s'<s$.
Let $0<\epsilon'<\epsilon$.  We  claim that
$$ \underset{0<|h|<\eta_0}\to\sup |D_{t,h}^{\epsilon'} u |_{\frac
12+\epsilon-\epsilon'}\le C|u|_{\frac 12
+\epsilon}.\tag 5.13$$ Assuming (5.13) for the moment, we get from (5.12)
that
$D_{t,h}^{\epsilon'} u \in W^{s'}$ for any $\frac 12 < s'<\frac
12+\epsilon-\epsilon' $.
Note that from (5.7), the operators $D_{t,h}^{\epsilon'}$ and
$\partial_b\db_b$ commute.
Thus we have
$$\partial_b\db_b D_{t,h}^{\epsilon'}u=D_{t,h}^{\epsilon'}\partial_b\db_b
u=D_{t,h}^{\epsilon'}\tilde\Theta^N\in C^{\alpha-\epsilon'}(U).\tag 5.14
$$
Applying the classical Schauder estimates to the equation (5.14) and
repeating the above
arguments used to obtain (5.11),
we have
$$  D_{t,h}^{\epsilon'}u \in L^\infty(U).\tag 5.15$$
This implies that $u\in C^{\epsilon'}(U)$.
To finish the proof of the lemma, it remains to prove  the claim (5.13). \qed
\enddemo
\proclaim{Lemma 5.3} Let $u$ be a function with $|u|_{s+\epsilon'}<\infty$
as above.
Then
$$\underset{0<|h|<\eta_0}\to \sup |D_{t,h}^{\epsilon'} u |_{s}\le
C|u|_{s+\epsilon'}. $$

\endproclaim
\demo{Proof}
We identify $h=(0,\cdots,0,h)$.
We have
$$ \aligned&\frac {\|D_{t,h}^{\epsilon'} u(
x+\tilde h )-D_{t,h}^{\epsilon'} u(x)\|}{|\tilde h |^s}\\&= \frac {\| ( u(
x +\tilde h +h )- u(x+\tilde h))-(u(x+h)-u(x))\|}{|h|^{\epsilon'}|\tilde h |^s}
\\&\le \frac {  (| u(
x+\tilde h   ) |_{s+\epsilon'}+ |u  |_{s+\epsilon'})|
h|^{s+\epsilon'}}{|h|^{\epsilon'} |\tilde h
|^s}
\\&\le \frac {  2|u  |_{s+\epsilon'} |  h|^{s+\epsilon'}}{|h|^{\epsilon'}
|\tilde h
|^s}.\endaligned\tag 5.16
$$
Similarly, we have
$$ \aligned&\frac {\|D_{t,h}^{\epsilon'} u(
x+\tilde h )-D_{t,h}^{\epsilon'} u(x)\|}{|\tilde h |^s}\\&= \frac {\| ( u(
x +\tilde h +h )-u(x+h))-( u(x+\tilde h)-u(x))\|}{|h|^{\epsilon'}|\tilde h |^s}
\\&\le \frac {  (| u(
x+  h   ) |_{s+\epsilon'}+ |u  |_{s+\epsilon'})|\tilde
h|^{s+\epsilon'}}{|h|^{\epsilon'} |\tilde h
|^s}
\\&\le \frac {  2|u  |_{s+\epsilon'} |  \tilde
h|^{s+\epsilon'}}{|h|^{\epsilon'} |\tilde h
|^s}.\endaligned\tag 5.17
$$

By  (5.16) and (5.17), we conclude that
$$\aligned \underset{0<|h|<\eta_0}\to\sup |D_{t,h}^{\epsilon'}
u|_{s}&=\underset{0<|h|<\eta_0}\to\sup\underset{0<|\tilde h
|<\eta_0}\to\sup \frac {\|D_{t,h}^{\epsilon'} u(
x+\tilde h )-D_{t,h}^{\epsilon'} u(x)\|}{|\tilde h |^s}
\\&\le C \underset{0<|h|+|\tilde h|<2\eta_0}\to\sup
\min(|h|^{s+\epsilon'},|\tilde h|^{s+\epsilon'}        ) \frac
{|u|_{s+\epsilon'}}{|h|^{s'}|\tilde h|^s}
\\&\le C|u|_{s+\epsilon'} .\endaligned
$$ This proves the lemma.  \qed
\enddemo
\demo{Proof of Theorem 1}
 Let $M$ be a  Levi-flat hypersurface of class $C^{2,\alpha}$ in $\Bbb
CP^n$, $n\ge 2$.  Let $u$ be obtained
 in Lemma 5.2 and
  $h= 2 \text{Im} u$, where $  \text{Im} u$ is the imaginary part of $u$.
From Lemma 5.2, the function
   $h$ is H\"older continuous on $M$.  Since $h$ is real-valued and $M$ is
compact,    there exists a point $Q_0$ such that
$h$ assumes a global
maximum at
$Q_0$. On the other hand, we have
$$\aligned \Theta_b & = \tilde{\Theta}^N|_{[T(M)]_{\Bbb R}
\cap J [T(M)]_{\Bbb R} }
\\ &= (-d\beta)|_{[T(M)]_{\Bbb R}
\cap J [T(M)]_{\Bbb R} }
 =-\partial_b\beta^{0,1}-\db_b\beta^{1,0}\\&=- \partial_b \db_b
u-\db_b\partial_b  \bar
u
\\&= \partial_b\db_b h\quad
 \endaligned\tag 5.18 $$
on $M$.
 From Proposition 1.2,  the curvature form   $i\Theta_b$ is a
positive $(1,1)$-form on $M$. Although $h $ is only  $C^{0, \epsilon}$-
continuous on $M$, its Hessian and $i\partial_b\db_b h$  can be computed by
using the barrier
functions, (e.g., see [Ca]).
  From
(5.18) and Proposition 1.2,
$h$ is a strictly plurisubharmonic function when restricted to each leaf
$\Sigma_t$  of
$M$. Hence,
at the global maximum point $Q_0 \in M $ of $h$, we obtain
a strictly plurisubharmonic function $h$  on that particular leaf $\Sigma_{0}$
containing $Q_0$ as an interior maximum. This
contradicts the maximum principle (cf. [Ca]). Thus there does not exist any
Levi-flat
hypersurface of class $C^{2,\alpha}$ in $\Bbb
CP^n$. Theorem 1 is proved. \qed

\enddemo

\enddocument